\documentclass[12pt]{article}
\usepackage{amssymb}
\usepackage[margin=1in]{geometry}
\usepackage{graphicx}
\usepackage{setspace}
\usepackage{svg}
\usepackage{mathtools}
\usepackage{amsmath}
\usepackage{xr}
\usepackage{caption}
\usepackage{todonotes}
\usepackage{caption}
\DeclareMathOperator*{\argmin}{arg\,min}
\usepackage{subcaption}
\usepackage{array}
 \usepackage{algorithm}
 \usepackage{algpseudocode}
\newcommand{\pg}[1]{{\color{black}  #1}}

\pdfoutput=1 
\title{Topological bifurcations in a mean-field game }
\author{Ali Akbar Rezaei Lori and  Piyush Grover\footnote{Email: Piyush.grover@unl.edu} \\ 
 Mechanical and Materials Engineering, \\University of Nebraska-Lincoln, Lincoln, Nebraska, USA}
 \date{\today}

\begin{document} 
\maketitle

\begin{abstract}
Mean-field games (MFG) provide a statistical physics inspired modeling framework for decision making in large-populations of strategic, non-cooperative agents. Mathematically, these systems consist of a forward-backward in time system of two coupled nonlinear partial differential equations (PDEs), namely the Fokker-Plank and the Hamilton-Jacobi-Bellman equations, governing the agent state and control distribution, respectively. In this work, we study a finite-time MFG with a rich global bifurcation structure using a reduced-order model (ROM). The ROM is a 4D two-point boundary value problem obtained by restricting the controlled dynamics to first two moments of the agent state distribution, i.e., the mean and the variance. Phase space analysis of the ROM reveals that the invariant manifolds of periodic orbits around the so-called `ergodic MFG equilibrium' play a crucial role in determining the bifurcation diagram, and impart a topological signature to various solution branches. We show a qualitative agreement of these results with numerical solutions of the full-order MFG PDE system.
\end{abstract}

\maketitle
\section{Introduction}

Mean-field game (MFG) theory \cite{bensoussan2013mean, caines2015mean} is a modeling framework for large-population non-cooperative engineering and socio-economic systems. This theory combines optimal control theory and game theory with ideas from statistical physics, where a `mean-field' approximation is used to simplify the study of systems with a large number of particles. In statistical physics the particles are generally passive, i.e., driven solely by environmental or inter-particle forces. On the other hand, MFGs are concerned with decision making or `active' particles (called agents). The MFGs are mathematically described by a coupled set of forward-backward in time nonlinear partial differential equations (PDEs).

The MFG PDEs have a fundamentally different structure than the PDEs used to describe most collective behavior in natural or engineering systems. These latter systems are generally modeled via reaction-diffusion or kinetic/hydrodynamics equations, etc., which are evolution equations solved forward in time. In the MFG systems, the collective behavior is the result of each agent solving  an optimal control problem that depends on its own state (and control), as well as the collective state \cite{bensoussan2013mean,nourian2013nash}. This imparts the MFG equations a forward-backward in time structure. 

In this work, we will work with a class of MFGs where the coupling between the agents is solely though the cost function. Consider a population of $N$ agents, where the vector-valued state $x_i\in\mathbb{R}^n$ of the $i$th agent is described by the following stochastic differential equation (SDE) :
 \begin{align}
dx_i(t)=f(x_i(t),a_i(t))dt+\sigma dW_i(t).\nonumber
\end{align}
Here $a_i\in\mathbb{A}$ is the control input, $\mathbb{A}$ is the set of permissible control inputs, and $W_i$ is the standard $n-$dimensional Wiener (noise) process. The combined state of the population can be described using the empirical distribution, whose density is $m(x,t)=\dfrac{1}{N}\sum\limits_{i=1}^{N}\delta(x-x_i(t))$. 
The control $a_i$ of each agent at each time is chosen to minimize the cost \begin{align}
    J(a_i)=\mathbb{E}\{\int_0^TG(x_i,a_i,x_{-i})dt+C(x_i(T))\big\},
\end{align}
where $G$ is the running cost that depends on the state and control of the $i$th agent as well as on the states of the other agents, and $C(x)$ is the terminal cost function. Here we use the notation $x_{-i}=\{x_1,x_2,\dots,x_{i-1},x_{i+1},\dots,x_N\}$. MFG theory postulates that the solution of the above problem in the limit $N\rightarrow\infty$ can be obtained by assuming an exogenous density $\hat{m}(x,t)$ during the optimization process of the $i$th agent \cite{bensoussan2013mean}. Hence, the running cost term is of the time-dependent form $\hat{G}(x_i,a_i;\hat{m}(.,t))$. This renders the optimization problem solvable by standard stochastic optimal control techniques \cite{borkar2005,yong2012stochastic}, and the optimal control is encoded in the value function $u(x_i,t)$ which solves the Hamilton-Jacobi-Bellman (HJB) equation:
\begin{align}
    \partial_tu(x,t)+\hat{H}(x,\nabla u(x,t),\hat{m}(.,t))+\frac{\sigma^2}{2}\Delta u(x,t)=0, \label{eq:HJB0}
\end{align}
with terminal condition $u(x,T)=C(x)$. Here, the Hamiltonian $\hat{H}(x,p,m)=\displaystyle\min_{a\in\mathbb{A}}\big(\hat{G}(x,a;\hat{m}(.,t))+p.f(x,a)\big)$, and the optimal feedback control can be obtained by $a(x,t)=\displaystyle\argmin_{a\in\mathbb{A}} \big(\hat{G}(x,a;\hat{m}(.,t))+\nabla u(x,t).f(x,a)\big)$ once $u(x,t)$ is known. The consistency condition of MFG requires that if each agent follows the above \emph{control law} $a(x,t)$, then the resulting agent density must equal the assumed exogenous density, i.e., $m(x,t)=\hat{m}(x,t)$. The density evolution is given by the Fokker-Plank (FP) equation
\begin{align}
\partial_tm(x,t)+\nabla(f(x,a(x,t))m(x,t))=\frac{\sigma^2}{2}\Delta m(x,t), \label{eq:FP0}
\end{align}
with initial condition $m(x,0)=m_0(x)$. To make the system well-posed, boundary conditions (in $x$) on $(u,m)$ also need to be specified. Equations (\ref{eq:HJB0},\ref{eq:FP0}), together with the respective final and initial conditions, and the boundary conditions, form the MFG PDE system. 
The MFG system consists of two coupled nonlinear PDEs with a forward-backward structure, i.e., the FP equation is well-posed forward in time, while the HJB is well-posed backward in time. A pair $(u,m)$ is the solution of the finite horizon MFG with time-horizon $T$ if it solves Eqs. (\ref{eq:HJB0},\ref{eq:FP0}), and satisfies the intial and final conditions as well as the boundary conditions.

Much of the research in the field of MFGs has focused on finding conditions on agent dynamics and cost functions to guarantee existence and uniqueness of solutions of these equations \cite{carmona2018probabilistic,bensoussan2013mean}. The uniqueness of solutions is generally guaranteed when the PDEs satsify certain monotonicity properties \cite{bardi2019non}. However, such monotonicity properties are satisfied only in special cases, and multiplicity of solutions is expected to be `generic'. 

This realization has led to recent interest in understanding and characterizing the multiple solution branches of non-monotonic MFG systems using bifurcation theory and abstract functional analytic techniques. Phase transitions in collective behavior such as flocking \cite{solon2015pattern}, synchronization of oscillators \cite{strogatz2000kuramoto}, and traffic flow 
\cite{gasser2004bifurcation}, are often studied via bifurcation analysis of evolution PDEs. In the same vein, bifurcations of the solutions of the MFG systems can be interpreted as phase transitions in the collective behavior of decision making agents.

In one of first works in this topic \cite{yin2012synchronization}, the authors applied bifurcation theoretic tools to a nonlocal infinite time-horizon MFG system posed on a periodic 1D domain, and proved the co-existence of a travelling wave time-periodic solution along with a steady (i.e., time-independent) MFG solution. Further rigorous analytic results for the same problem were recently obtained in \cite{carmona2023synchronization,cesaroni2024stationary}. In \cite{grover2018mean}, a bifurcation leading to multiple co-existing steady MFG solutions was demonstrated in an infinite time-horizon MFG model posed on the whole real line. In certain classes of non-monotonic finite time-horizon MFGs, abstract bifurcation theory \cite{kielhofer2006bifurcation} along with eigenfunction expansions in periodic spatial domains were applied to prove the existence of time-periodic solutions \cite{cirant2019existence,cirant2018variational}. Another recent work uses Aubry-Mather theory for proving the existence of time-periodic solutions in first-order MFGs \cite{ni2024time}.

In addition to analytic techniques employed in the above mentioned works, the use of exactly solvable and reduced-order models \cite{hongler2020mean,swiecicki2016schrodinger, ullmo2019quadratic} along with geometric tools of dynamical systems theory \cite{wiggins2003introduction} has provided valuable insight into solution regimes of non-monotonic MFGs. In \cite{swiecicki2016schrodinger, ullmo2019quadratic}, an exact mapping of the so-called quadratic MFGs into a pair of Schrodinger equations was exploited to derive low-dimensional two-point boundary value problem (BVP) models under various assumptions on relative strengths of different terms in the MFG system, as well as on initial density. For instance, if the initial agent density is Gaussian, and the attractive interaction term is dominant, the solution is expected to stay approximately Gaussian for all times. In this case, a 4D BVP model was derived in \cite{ullmo2019quadratic} where the state vector consists of the population mean and standard deviation, and their corresponding momentum variables. A key advantage of such low dimensional models is that their behaviour can be understood using geometric tools of phase space analysis. In \cite{ullmo2019quadratic}, by assuming that the mean and variance dynamics are decoupled, the invariant manifolds of the equilibrium point of a 2D reduced-order model of the infinite-time horizon (`ergodic') problem were used to obtain qualitative and quantitative results for solutions of the corresponding full-order finite-time horizon MFG PDE system.

The goal of this paper is to derive and analyze a related reduced-order 4D Hamiltonian two-point BVP model of a finite-horizon MFG system that exhibits multiplicity of solutions for large enough time-horizon $T$, and relate the origin of various solution branches to the 4D phase space geometry. We demonstrate that the different solution branches are topologically distinct, and that these solutions can be understood as transitions through a bottleneck in the configuration space between the initial and final conditions. The phase space flow in the bottleneck is organized by the cylindrical 2D stable and unstable manifolds of an unstable periodic orbit that exists in the bottleneck. Similar transitions dynamics have been previously studied in the context of the gravitational three-body problem \cite{mcgehee1969some,koon2000heteroclinic}, chemical kinetics \cite{de1991cylindrical}, and structural mechanics \cite{zhong2018tube}. The qualitative (i.e., topological) aspects of the results from analysis of the BVP are shown to persist in numerically computed solutions of the full-order MFG PDE.

\section{The Quadratic MFG PDE and Reduced Order Modeling}

We begin by recalling the setup of \cite{ullmo2019quadratic}, and consider the case where the scalar state $x_i\in\mathbb{R}$ of $i$th agent is driven by a control term and standard Brownian noise as :
\begin{equation}\
    dx_i(t)=a_i(t)dt+\sigma dw_i(t),
    \label{SDE}
\end{equation}
and the running cost $G=R-\bar{V}$ consists of two parts: a quadratic control penalization term $R(x_i,a_i)=\dfrac{\mu a_i^2}{2}$, and a potential cost $\bar{V}[m](x_i)$ that depends on density of agents $m$. The potential cost is further split into interaction and external potentials as $\bar{V}[m](x_i)=f[m](x_i)+U_0(x_i)$, where 
\begin{equation}
    \begin{split}
        \label{costs}
        & f[m](x_i)=g\,m^\alpha(x_i) ,\\&
        U_0(x)=-h\frac{x_i^2}{2}-\frac{x_i^4}{4}.
    \end{split}
\end{equation}

The resulting HJB and FP equations are:
\begin{align}
    &\partial_tu(x,t)-\frac{1}{2\mu}(\partial_xu(x,t))^2+\frac{\sigma^2}{2}\partial_{xx}u(x,t)=\bar{V}[m](x,t),\label{eq:HJB1}
    \\ & \partial_tm(x,t)-\frac{1}{\mu}\partial_x(m(x,t)\partial_xu(x,t))-\frac{\sigma^2}{2}\partial_{xx}m(x,t)=0,\label{eq:FP1}
\end{align}
with a prescribed initial density $m(x,0)=m_0(x)$, and a prescribed terminal value function $u(x,T)=C(x(T))$. The corresponding feedback control \pg{$a_i(t)=-\dfrac{1}{\mu}\partial_xu(x_i,t)$}.

\paragraph{Ergodic MFG:} In the limit of infinite time horizon ($T\to\infty$), the MFG equations reduce to the following steady state equations:
\begin{equation}
    \begin{split}
        &-\lambda^e-\frac{1}{2\mu}(\partial_xu^e(x))^2+\frac{\sigma^2}{2}\partial_{xx}u^e(x)=\bar{V}[m^e](x),
    \\ & \frac{1}{\mu}\partial_x(m^e(x)\partial_xu^e(x))+\frac{\sigma^2}{2}\partial_{xx}m^e(x)=0,\label{eq:MFG_ergodic}
    \end{split}
\end{equation}
where the pair $ (m^e(x), u^e(x))$  is the ergodic equilibrium, and $\lambda^e$ is the ergodic constant. It has been shown  \cite{cardaliaguet2013long} that solution of the long time horizon (i.e., $T\gg 1$) finite time MFG (Eqs. \ref{eq:HJB1}, \ref{eq:FP1}) tends to approach and stay close to the solution $(m^e(x), u^e(x))$ during the interval time $0\ll t\ll T$. This property is a recurring theme in calculus of variations and optimal control problems, and is commonly referred to as the `turnpike property' \cite{zaslavski2005turnpike}.
 
In \cite{ullmo2019quadratic}, Cole-Hopf transformations $\Phi(x,t)=\exp(\dfrac{-u(x,t)}{\mu\sigma^2})$ and $\Gamma(x,t)=\dfrac{m(x,t)}{\Phi(x,t)}$, were used to reduce the finite-time MFG system (Eqs. \ref{eq:HJB1}, \ref{eq:FP1}) into the following pair of nonlinear diffusions:

\begin{align}
-\mu\sigma^2\partial_t\Phi(x,t)=\frac{\mu\sigma^4}{2}\partial_{xx}\Phi(x,t)+\bar{V}[m](x,t)\Phi(x,t), \label{NLS1}  \\ 
\mu\sigma^2\partial_t\Gamma(x,t)=\frac{\mu\sigma^4}{2}\partial_{xx}\Gamma(x,t)+\bar{V}[m](x,t)\Gamma(x,t). \label{NLS2}
\end{align}
By exploiting the analogy of the above equations with imaginary time nonlinear Schrodinger equation, the following expressions for solutions of the finite-time MFG system in terms of the ergodic solution pair $(m^e(x), u^e(x))$ and ergodic constant $\lambda^e$ were obtained: 
\begin{align}
    & \Phi(x,t)=\exp\big({\frac{\lambda^et}{\mu\sigma^2}}\big)\psi^e(x), \label{eq:phi_gamma_sol1} \\
    & \Gamma(x,t)=\exp\big({-\frac{\lambda^et}{\mu\sigma^2}}\big)\psi^e(x).  \label{eq:phi_gamma_sol2}
\end{align}

\paragraph{Variational Principle}

Furthermore, it was shown in \cite{ullmo2019quadratic} that the Eqs.(\ref{NLS1},\ref{NLS2}) can be obtained as necessary conditions for the stationarity w.r.t $(\Phi,\Gamma)$ of the following action functional :
\begin{equation}
    \label{action_fncl}\begin{split}
    \mathcal{S}[\Phi,\Gamma]& =\int_0^t\int_\mathbb{R}\Bigg[-\frac{\mu\sigma^2}{2}\big((\partial_t\Gamma)\Phi-(\partial_t\Phi)\Gamma\big)-\frac{\mu\sigma^4}{2}(\partial_x\Phi)(\partial_x\Gamma)+\Phi U_0(x)\Gamma+F[\Phi\Gamma]\Bigg] dx dt
        \end{split} = \int_{0}^t L dt, 
\end{equation}
where $F[m]=\displaystyle\int\limits^m f(m')\,dm' $. Here the Lagrangian $L=\displaystyle\int\limits_\mathbb{R}\big[-\dfrac{\mu\sigma^2}{2}\big(\partial_t(\Gamma)\Phi-\partial_t(\Phi)\Gamma\big)\big]dx+E_{tot}$. The total energy $E_{tot}=E_{kin}+E_{ipot}+E_{epot}$, is the sum of the kinetic term $E_{kin}=\dfrac{\mu\sigma^4}{{2}}\displaystyle\int_{\mathbb{R}}\Phi\partial_{xx}\Gamma dx$ (using integration by parts), the interaction potential $E_{ipot}=\displaystyle\int_{\mathbb{R}} F[m]dx$ and the external potential $E_{epot}=\displaystyle\int_{\mathbb{R}} U_0(x)m(x,t)dx$.

\subsection{Derivation of the BVP system in Lagrangian variables}\label{subsec:deriv_BVP}

We consider the scenario where the interaction potential is attractive, and dominates the external potential. In this case,
a Gaussian distribution with density $m^G(x,t) = \dfrac{1}{\sqrt{2\pi \epsilon^2S^2}}\exp\left(-\dfrac{(x-X)^2}{2\epsilon^2S^2}\right)$ is a good approximation for a population of agents with mean $X(t)=\displaystyle\int x\,m(x,t)\,dx$, and standard deviation $\Sigma(t)=\epsilon S(t)=\sqrt{\displaystyle\int x^2\,m(x,t)\,dx-X^2(t)}$. Here, we have introduced the parameter $\epsilon$ (with $0<\epsilon<1$), and use scaled standard deviation $S$ in order to facilitate the search for appropriate system parameters in the following sections. In the spirit of Eqs. (\ref{eq:phi_gamma_sol1}, \ref{eq:phi_gamma_sol2}), and following \cite{pitaevskii2016bose,ullmo2019quadratic}, we use an anstaz for $\Phi(x,t)$ and $\Gamma(x,t)$ shown below:
\begin{equation}
    \begin{split}\label{ansatz}
        \Phi(x,t)=\exp{\big(\frac{-\gamma+Px}{\mu\sigma^2}\big)}\frac{1}{(2\pi\epsilon^2S^2)^\frac{1}{4}}\exp{\big(-\frac{(x-X)^2}{4\epsilon^2S^2}(1-\frac{\Lambda}{\mu\sigma^2})\big)},\\
        \Gamma(x,t)=\exp{\big(\frac{\gamma-Px}{\mu\sigma^2}\big)}\frac{1}{(2\pi\epsilon^2S^2)^\frac{1}{4}}\exp{\big(-\frac{(x-X)^2}{4\epsilon^2S^2}(1+\frac{\Lambda}{\mu\sigma^2})\big)},
    \end{split}
\end{equation}
 where $\gamma(t)=\dfrac{\Lambda(t)}{4}+\gamma_0$ is a necessary condition for the above ansatz to satisfy Eqs. (\ref{NLS1},\ref{NLS2}). Hence, $\Phi(x,t)$ and $\Gamma(x,t)$ are effectively parameterized using four time-varying scalar variables: $Z=[X,S,P,\Lambda]^\intercal$. Given a density and value function pair $(m,u)$, we can compute $P(t)=\displaystyle\int -\mu\sigma^2\Phi(x,t)\,\partial_x(\Gamma(x,t))\,dx$, and $\Lambda(t)=\displaystyle\int -\mu\sigma^2\Phi(x,t)\bigg(x\,\partial_x\Gamma(x,t)+\partial_x\big(x\Gamma(x,t)\big)\bigg)\,dx-2XP$, where we recall that $\phi(x,t)$ and $\Gamma(x,t)$ are obtained via the Cole-Hopf transformation.

The kinetic and interaction potential energies in this case are:
\begin{equation}\label{energies2}
    \begin{split}
        &E_{kin}=\frac{{P }^2 }{2\mu}+\frac{{{\Lambda }^2 }} {8\mu\,\epsilon^2S^2 }-\frac{{\mu\,\sigma^4 }}{8\,\epsilon^2S^2 }
        ,\\& E_{ipot}=\frac{g}{(2\pi)^{\frac{\alpha}{2}} {{\left(\alpha+1\right)}}^{3/2}\,(\epsilon\,S)^\alpha}\,.
    \end{split}
\end{equation}

We approximate the external potential energy using a Taylor expansion of $U_0(x)$ around the mean $X$ as follows:\pg{
\begin{equation}\label{potential expansion}\begin{split}
    &\int U_0(x)\,m^G(x,t)\,dx=\, U_0(X)+\int (x-X)\partial_xU_0(X)\,m^G(x,t)\,dx +\frac{1}{2}\int(x-X)^2\,\partial_x^2U_0(X)m^G(x,t)\,dx\\&+\frac{1}{3!}\int \big(x-X)^3\partial_x^3U_0(X)\,m^G(x,t)\,dx +\frac{1}{4!}\int \big(x-X)^4\partial_x^4U_0(X)\,m^G(x,t)\,dx +\dots .  
\end{split}
\end{equation} }
Using the fact that for the Gaussian distribution $m^G$, 
\begin{equation}
    \label{Gau_Moment}
    \int(x-X)^nm^G\,dx =\begin{cases} (n-1)!!(\epsilon S)^n & \text{for n even}, \\
                     0 &  \text{for n odd},
       \end{cases}
\end{equation}
where $(n-1)!!$ denotes the product of all odd numbers less than $n$, and keeping terms up to the fourth order, we get:
\begin{equation}
    \label{E_pot}
    E_{epot}=U_0(X)+\frac{1}{2}(\epsilon S)^2\partial^2_xU_0(X)+\frac{1}{4!}\partial^4_xU_0(X)\,3(\epsilon S)^4.
\end{equation}

By substituting Eqs.(\ref{energies2},\ref{E_pot}) and the ansatz Eq.(\ref{ansatz}) in Eq.(\ref{action_fncl}), we obtain the reduced functional $\Bar{\mathcal{S}}=\displaystyle\int \bar{L}(X,\dot{X},P,S,\dot{S},\Lambda) \,dt$. The reduced Lagrangian $\bar{L}$ is:
\begin{equation}\label{Lagrang_rescale}
    \begin{split}
        &\bar{L}(X,\dot{X},P,\Sigma,\dot{\Sigma},\Lambda)=-P\dot X-\frac{\Lambda \dot S }{2 S }+U_0(X)+\frac{1}{2}(\epsilon S)^2\partial^2_xU_0(X)+\frac{1}{4!}\partial^4_xU_0(X)\,3(\epsilon S)^4\\&+\frac{{P }^2 }{2\mu}+\frac{{{\Lambda }^2 -\sigma^4 \,\mu^2 }}{8\mu\,{(\epsilon S) }^2 }+\frac{g}{(2\pi)^{\frac{\alpha}{2}} {{\left(\alpha+1\right)}}^{3/2}\,{(\epsilon S) }^\alpha},
    \end{split}
\end{equation} 
where we have used the relation $\dot\gamma=\dfrac{\dot\Lambda}{4}$.

Using Euler-Lagrange equations $\dfrac{d}{dt}(\dfrac{\partial L}{\partial \dot{Z}_i})-\dfrac{\partial L}{\partial Z_i}=0$ for $i=1,2,3,4$, we obtain the following coupled system of equations:
\begin{equation}\label{higher_sigma1}
    \begin{split}
         &\dot{X}=\frac{P}{\mu},\\&
    \dot{P}=X^3+hX+3(\epsilon S)^2X ,\\&\dot{S}=\frac{\Lambda}{2\mu\epsilon^2 S},\\&\Dot{\Lambda}=\frac{\Lambda^2-\mu^2\sigma^4}{2\mu(\epsilon S)^2}+\frac{2g\alpha}{\alpha+1}\frac{1}{\sqrt{\alpha+1}(2\pi)^{\alpha/2}}\frac{1}{(\epsilon S)^\alpha}+2(\epsilon S)^2(3X^2+h)+6(\epsilon S)^4.
    \end{split}
\end{equation}

The first two equations for $X$ and $P$ govern the evolution of the mean, and the last two equations for $S$ and $\Lambda$ govern the variance of the Gaussian distribution. The BVP system is obtained by appending the above system of equations with the initial and final time conditions as follows: 
\begin{equation}\label{IC and FC1}
    \begin{split}
       &X(0)=\int x\,m(x,0)\,dx,\\&
       X(T)=\int x\,m(x,T)\,dx,
       \\& S(0)=\frac{1}{\epsilon}\sqrt{\int x^2\,m(x,0)\,dx-(X(0))^2},\\&
      S(T)=\frac{1}{\epsilon}\sqrt{\int x^2\,m(x,T)\,dx-(X(T))^2}.
    \end{split}
\end{equation}
\subsection{Transformation of the BVP system into Hamiltonian variables}
\label{subsec:deriv_ham}
To facilitate the phase space analysis of the BVP system, we will work with Hamiltonian rather than Lagrangian variables. We use the Legendre transformation using configuration variables $(q_1=X,q_2=S)$ to obtain the corresponding conjugate momenta:
\begin{equation}\label{L to H}
    \begin{split}
     p_1=\frac{\partial L}{\partial \dot{q_1}}=-P,\\ p_2=\frac{\partial L}{\partial \dot{q_2}}=-\frac{\Lambda}{2q_2}.
    \end{split}
\end{equation}

The Hamiltonian $H=-P\dot X-\frac{\Lambda\dot S}{2S}- L=-E$ is: 
\begin{equation}\begin{split}\label{ham_higherS4}
    H(q_1,p_1,q_2,p_2)=-\frac{p_1^2}{2\mu}-\frac{p_2^2}{2\epsilon^2\mu}-V(q_1,q_2),
    \end{split}
\end{equation}
where
\begin{equation}\label{eq:pot_couple}
    \begin{split}
       V(q_1,q_2)=-h\frac{q_1^2}{2}-\frac{q_1^4}{4}-\frac{1}{2}\,(\epsilon q_2)^2\,(h+3q_1^2)-\frac{\mu\sigma^4}{8 (\epsilon q_2)^2}+\frac{g}{(\alpha+1)\sqrt{\alpha+1}\,(2\pi)^{\alpha/2}\,(\epsilon q_2)^\alpha}-\frac{3}{4}(\epsilon q_2)^4
    \end{split}
\end{equation}
is the potential energy.

The resulting Hamiltonian equations of motion are:
\begin{equation}\label{eq:Ham_ode}
    \begin{split}
         &\dot{q_1}=\frac{\partial H}{\partial p_1}=-\frac{p_1}{\mu},\\&
    \dot{p_1}=-\frac{\partial H}{\partial q_1}=-q_1^3-h\,q_1-3\,(\epsilon q_2)^2\,q_1 ,\\&\dot{q_2}=\frac{\partial H}{\partial p_2}=-\frac{p_2}{\epsilon^2\mu},\\&\Dot{p_2}=-\frac{\partial H}{\partial q_2}=\frac{\mu\epsilon\sigma^4}{4(\epsilon q_2)^3}-\frac{g\epsilon\alpha}{\alpha+1}\frac{1}{\sqrt{\alpha+1}(2\pi)^{\alpha/2}}\frac{1}{(\epsilon q_2)^{\alpha+1}}-\epsilon^2\, q_2(3\,q_1^2+h)-3\epsilon^4 q_2^3.
    \end{split}
\end{equation}
The Hamiltonian BVP system is obtained by appending the above system with the following initial and final time conditions:
\begin{equation}\label{eq:ic_fc_ham}
    \begin{split}
       &q_1(0)=X(0),\\&
       q_1(T)=X(T),
       \\& q_2(0)=S(0),\\&
      q_2(T)=S(T),
    \end{split}
\end{equation}
where the r.h.s. of the above equations are specified in Eqs. \ref{IC and FC1}.

\section{Phase space geometry of the Hamiltonian ODE}
Since the Hamiltonian BVP system defined in Eqs. (\ref{eq:Ham_ode},\ref{eq:ic_fc_ham}) is a reduced-order model of the finite-horizon MFG, the equilibria of the Hamiltonian system of Eqs. \ref{eq:ic_fc_ham} correspond to the ergodic equilibria, i.e., the solutions of infinite time MFG described by Eqs. \ref{eq:MFG_ergodic}. As a consequence of the turnpike property, the solutions of BVP with $T\gg 1$ spend most of their time near the equilibria of Eqs. \ref{eq:ic_fc_ham}, except at the beginning and the end of the time-horizon. In this section, we analyze the phase space geometry around equilibria of the Hamiltonian ODE.

The linear stability of an equilibrium point $\mathcal{X}=\begin{pmatrix}q_1&0&q_2&0\end{pmatrix}^\intercal$ of Eqs. \ref{eq:ic_fc_ham} is determined by the Jacobian matrix:

\begin{equation}\hspace{-0.4in}
    J_h|_{\mathcal{X}}=\begin{array}{l}
\left(\begin{array}{cccc} 0 & -\dfrac{1}{\mu } & 0 & 0\\ -3{\epsilon }^2 \,{q_2 }^2-3\,{q_1 }^2 -h & 0 & -6\,{\epsilon }^2 \,q_1 \,q_2  & 0\\ 0 & 0 & 0 & -\dfrac{1}{\mu\epsilon^2 }\\-6\,\epsilon^2 \,q_1 \,q_2  & 0 & \dfrac{\alpha\,\epsilon^2\,g}{{{\left(\epsilon \,q_2 \right)}}^{\alpha+2} \,{{\left(2\,\pi \right)}}^{\alpha/2} \,\sqrt{\alpha+1}}-9\,{\epsilon }^4 \,{q_2 }^2 -\dfrac{3\,\sigma^4 \,\mu }{4\,{\epsilon }^2 \,{q_2 }^4 }-\epsilon^2 \,{\left(3\,{q_1 }^2  +h\right)} & 0
\end{array}\right).
\end{array}\label{eq:jac}
\end{equation}

\begin{figure}[hbt!]
\includegraphics[width=0.49\textwidth]{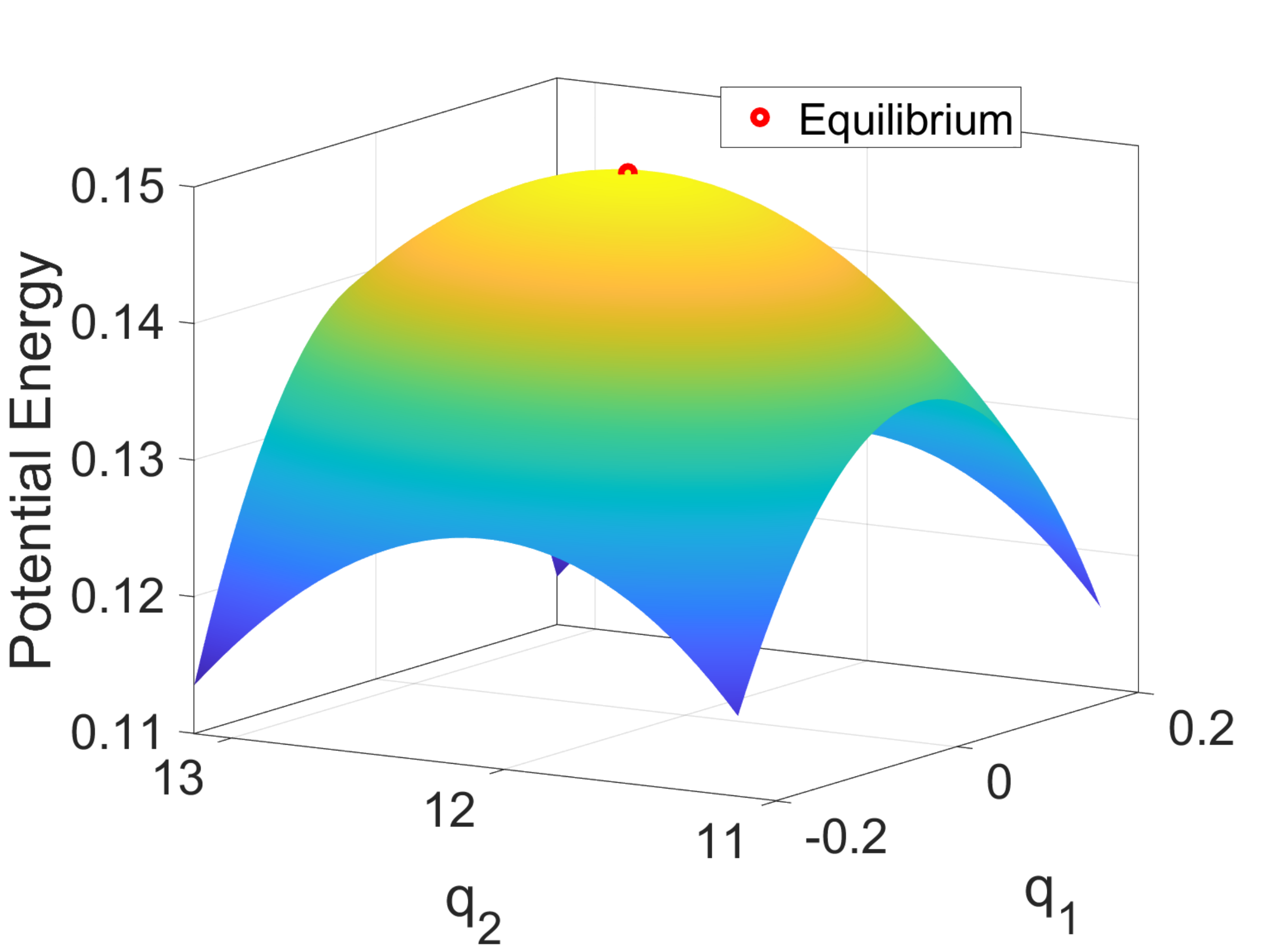} 
\includegraphics[width=0.49\textwidth]{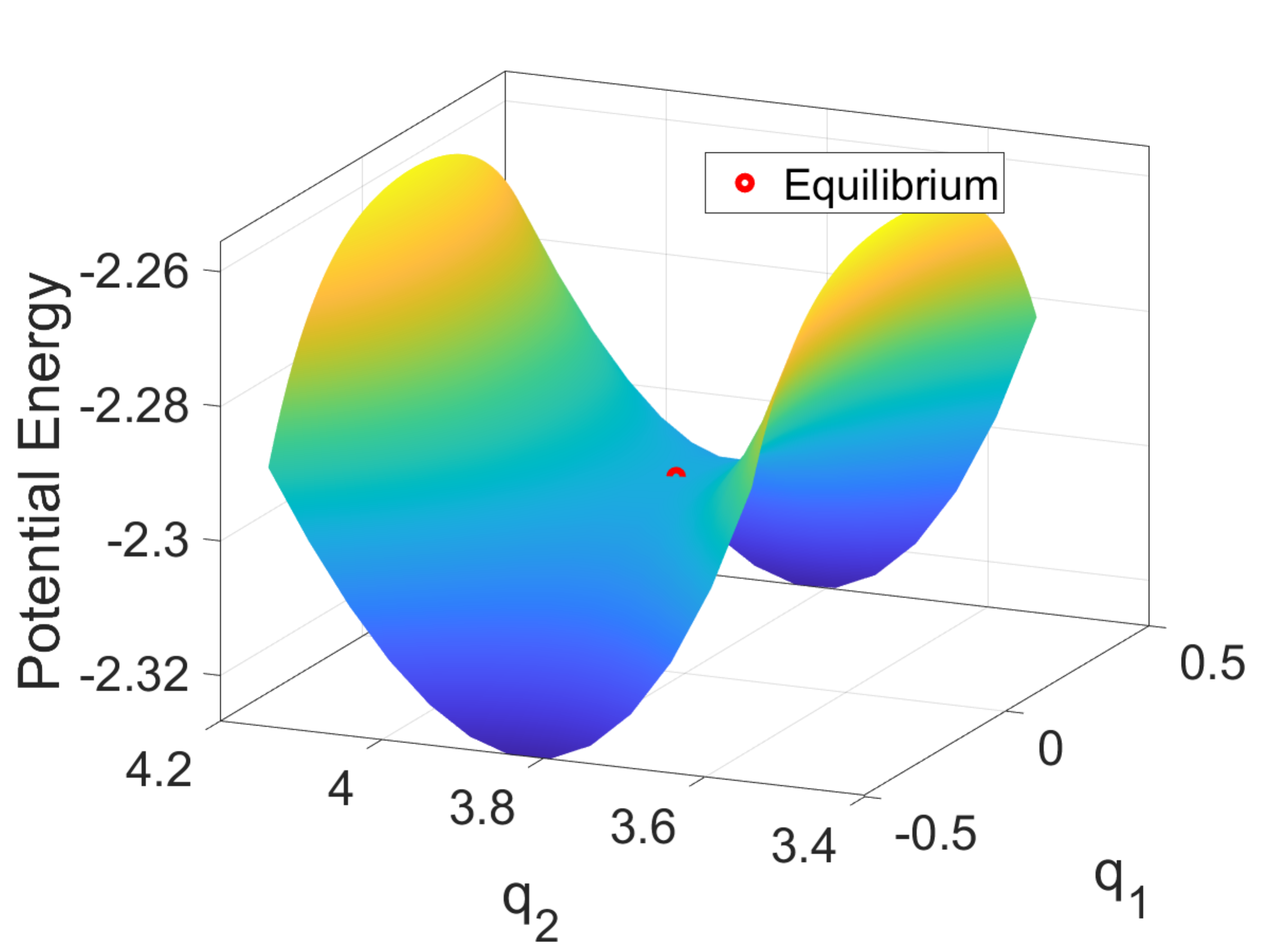} 
\caption{\footnotesize Potential energy surface in (Left) $saddle\times saddle$ and (Right) $saddle\times center$ case. }
\label{fig:pot_both}
\end{figure}

We consider two types of equilibria, see Fig. \ref{fig:pot_both}:    
\begin{enumerate}
    \item $Saddle\times saddle$ type: We call an equilibrium $\mathcal{X}=\begin{pmatrix}q_1&0&q_2&0\end{pmatrix}^\intercal$ to be of $saddle\times saddle$ type if eigenvalues of $J_h|_{\mathcal{X}}$ are of the form $\lambda_{1,2}=\pm \gamma_1, \lambda_{3,4}=\pm \gamma_2$, where $\gamma_1$ and $\gamma_2$ are both real. In this case, the pair $(q_1,q_2)$ is a local maxima of the potential energy function $V(q_1,q_2)$ defined in Eq. \ref{eq:pot_couple}.
    \item $Saddle\times center$ type: We call an equilibrium $\mathcal{X}=\begin{pmatrix}q_1&0&q_2&0\end{pmatrix}^\intercal$ to be of $saddle\times center$ type if eigenvalues of $J_h|_{\mathcal{X}}$ are of the form $\lambda_{1,2}=\pm \lambda, \lambda_{3,4}=\pm i\nu$, where $\lambda$ and $\nu$ are both real. In this case, the pair $(q_1,q_2)$ is a saddle point of the potential energy function $V(q_1,q_2)$ defined in Eq. \ref{eq:pot_couple}.
    
\end{enumerate}

\subsection{Phase space geometry near a $saddle\times saddle$ type ergodic equilibrium}

The phase space geometry near a $saddle\times saddle$ type equilibrium is governed by the stable and unstable eigenvectors of the Jacobian matrix $A\triangleq J_h|_{\mathcal{X}}$:
\begin{equation}
    A=\begin{pmatrix}
0 & -a & 0 & 0\\
-b & 0 & 0 & 0\\
0 & 0 & 0 & -c \\
0 & 0 & -d & 0
\end{pmatrix},\label{eq: jac SS}
\end{equation}
where $a,b,c,d$ are all positive real numbers.
Suppose $(u_1,v_1)$ and $(u_2,v_2)$ are the eigenvectors of $A$ corresponding to eigenvalues $\lambda_{1,2}=\pm\gamma_1\triangleq\pm\sqrt{ab}$ and $\lambda_{3,4}=\pm\gamma_2\triangleq\pm\sqrt{cd}$, respectively. Let $\mathcal{Z}=\begin{pmatrix} q_{1l}&p_{1l}&q_{2l}&p_{2l}\end{pmatrix}^\intercal$ denote perturbation about $\mathcal{X}$. Then, the linearized system $\dot{\mathcal{Z}}=A\mathcal{Z}$ can be written in the eigenvector basis $(u_1,v_1, u_2,v_2)$ as follows:
\begin{equation}
    \dot{\mathcal{Y}}=D\mathcal{Y}=\begin{pmatrix}
\gamma_1 & 0 & 0 & 0\\
0 & -\gamma_1 & 0 & 0\\
0 & 0 & \gamma_2 & 0 \\
0 & 0 & 0 & -\gamma_2
\end{pmatrix}\mathcal{Y},
\end{equation}
where $\mathcal{Z}=T\mathcal{Y}$ and $D=T^{-1}AT$, with the eigenvector matrix
\begin{equation}
T=\begin{pmatrix}
    \sqrt{\dfrac{a}{a+b}}& \sqrt{\dfrac{a}{a+b}}&0&0\\
     -\sqrt{\dfrac{b}{a+b}} &\sqrt{\dfrac{b}{a+b}}&0&0\\
     0 &0& \sqrt{\dfrac{c}{c+d}} &\sqrt{\dfrac{c}{c+d}}\\
     0&0&-\sqrt{\dfrac{d}{c+d}}&\sqrt{\dfrac{d}{c+d}}
\end{pmatrix}.\label{eq:T_ss}
\end{equation}

\begin{figure}[hbt!]
\includegraphics[width=.49\textwidth]{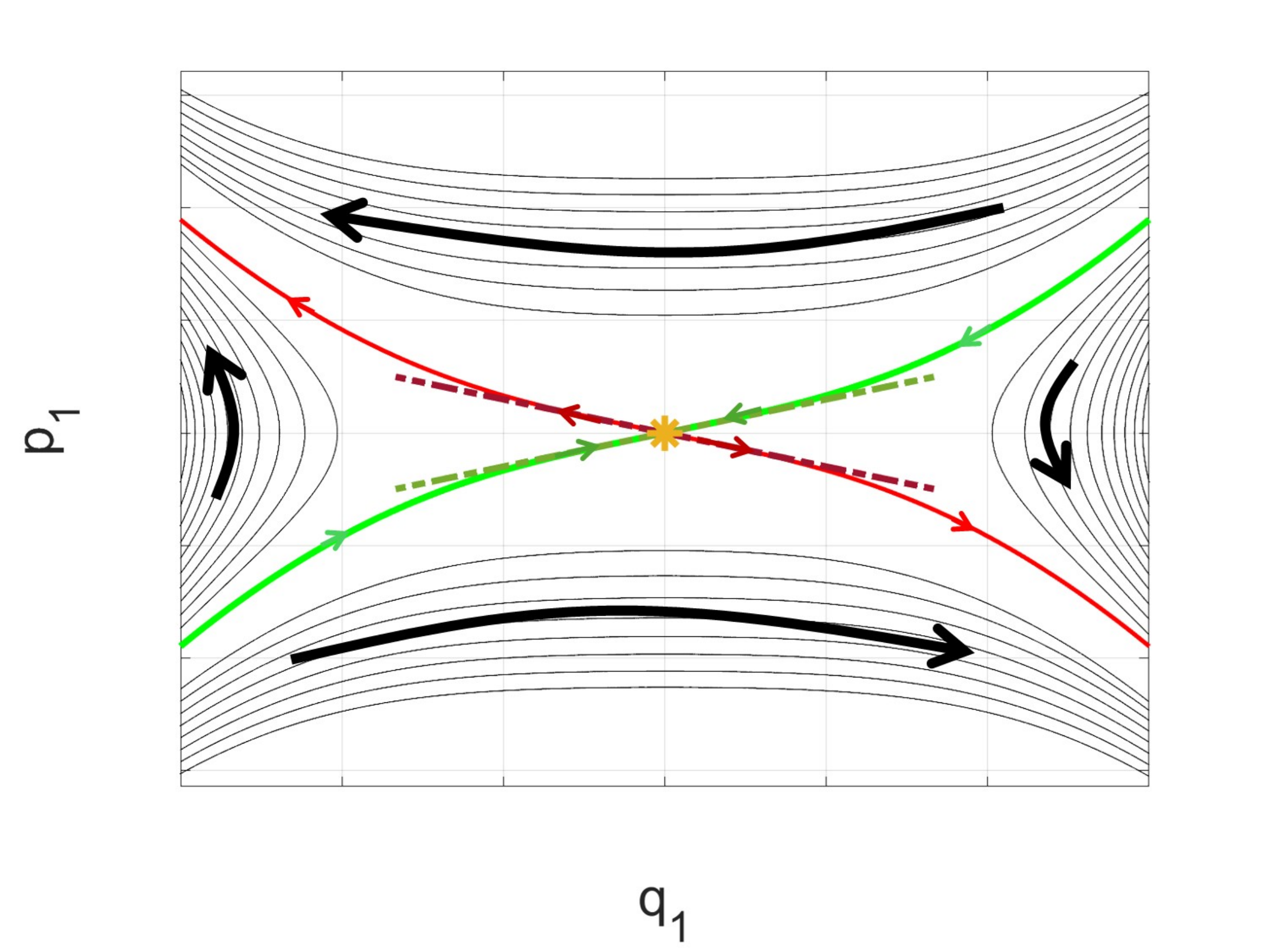}
\includegraphics[width=.49\textwidth]{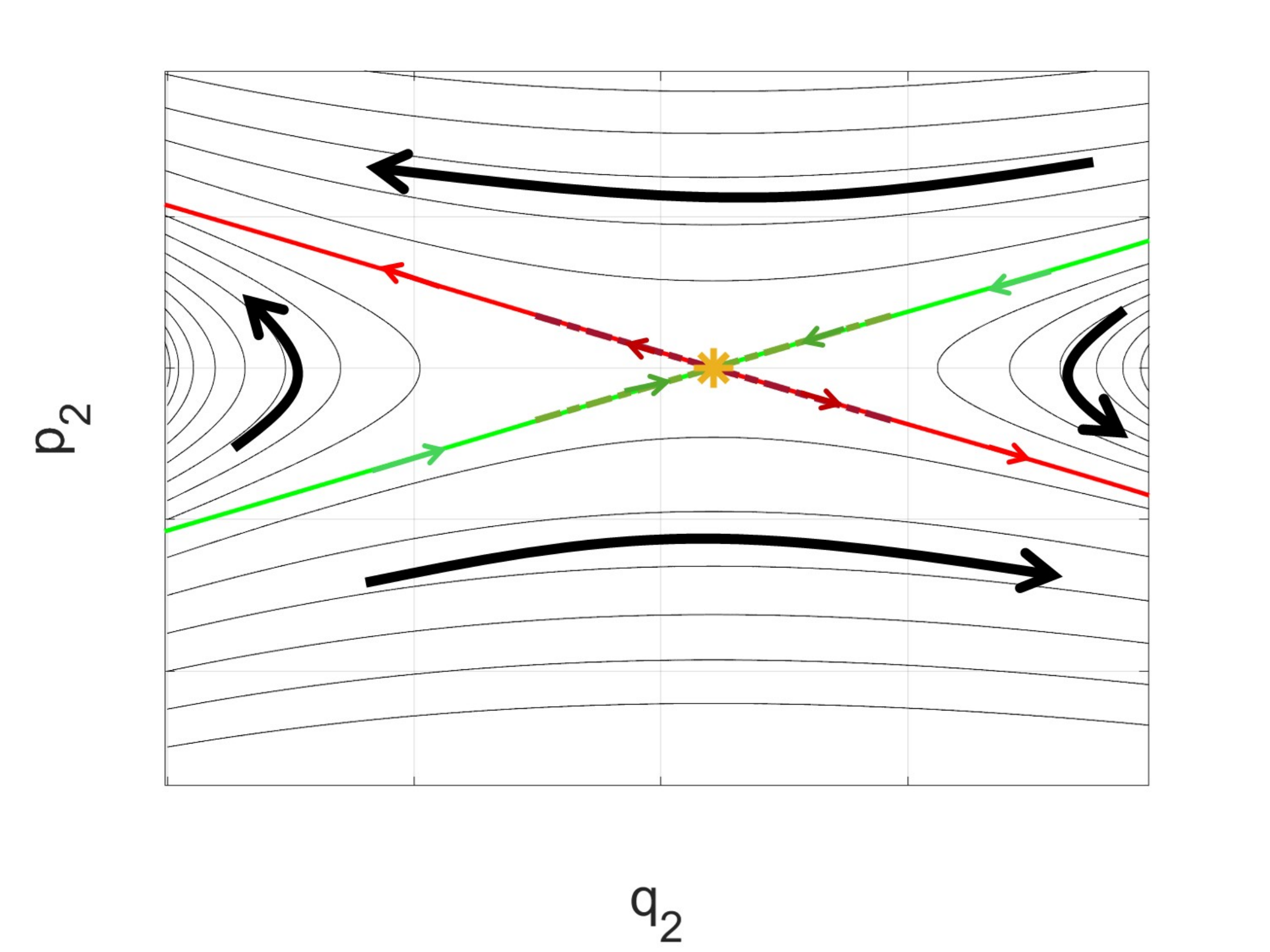} 
\caption{\footnotesize{Phase portrait of the nonlinear Hamiltonian ODEs in the neighborhood of a $saddle\times saddle$ equilibrium, projected on the  (Left) $q_1-p_1$ plane, and (Right) $q_2-p_2$ plane. Also shown are stable (green) and unstable (red) eigenvectors, and the corresponding invariant manifolds.}}
\label{fig:S-S_ODE}
\end{figure}

Hence, linearization of the Hamiltonian system (Eqs.\ref{eq:Ham_ode}) yields two decoupled systems in the two eigenspaces. In the nonlinear system, the flow is governed by the two dimensional stable manifold $W^s$, and the two-dimensional unstable invariant manifold $W^u$. Here, 
$W^s=\{Z\in\mathbb{R}^4|\lim_{t\rightarrow\infty}\phi_t(Z)\rightarrow \mathcal{X}\}$, and $W^u=\{Z\in\mathbb{R}^4|\lim_{t\rightarrow -\infty}\phi_t(Z)\rightarrow \mathcal{X}\}$, and $\phi_t$ is the time$-t$ flowmap of Eqs. \ref{eq:Ham_ode}. The stable (resp. unstable) manifold is tangent to the stable eigenspace (resp. unstable eigenspace).  Fig. \ref{fig:S-S_ODE} shows the phase space near the equilibrium projected on to the $q_1-p_1$ and $q_2-p_2$ planes.

\subsection{Phase space geometry near a $saddle\times center$ type ergodic equilibrium}
\label{subsec:phase_sc}

\begin{figure}[hbt!]
\includegraphics[width=.3\textwidth]{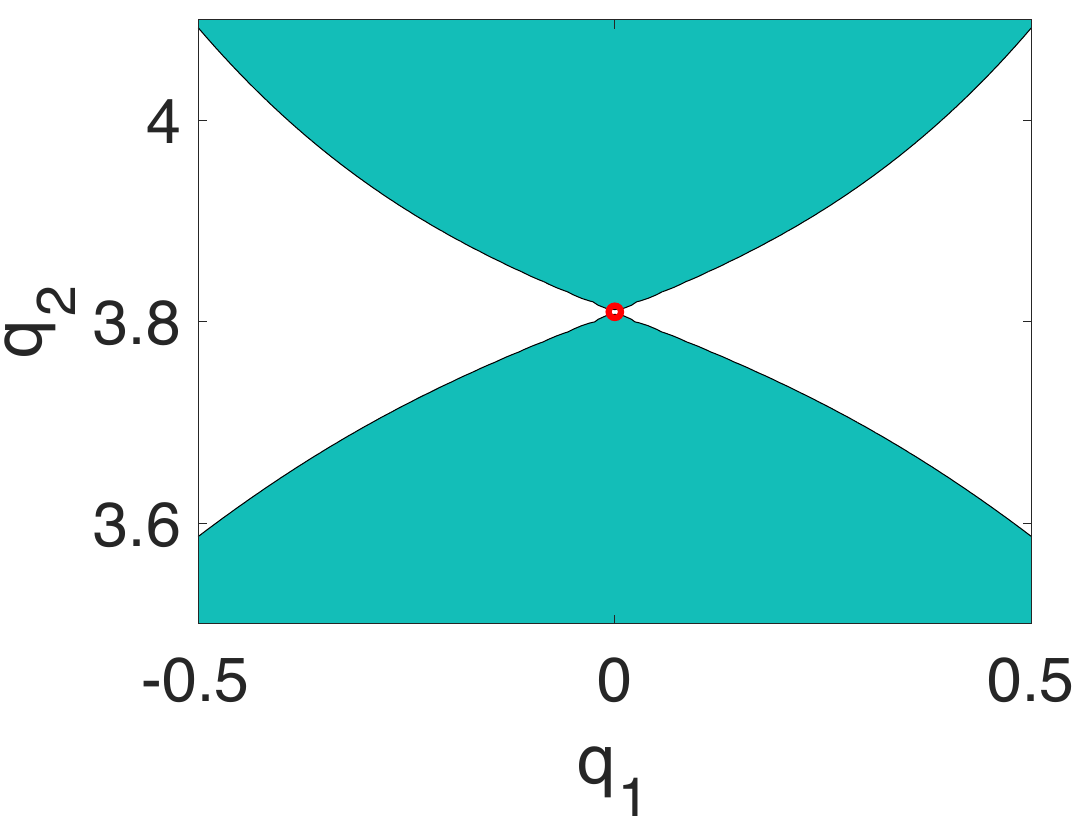}
\includegraphics[width=.3\textwidth]{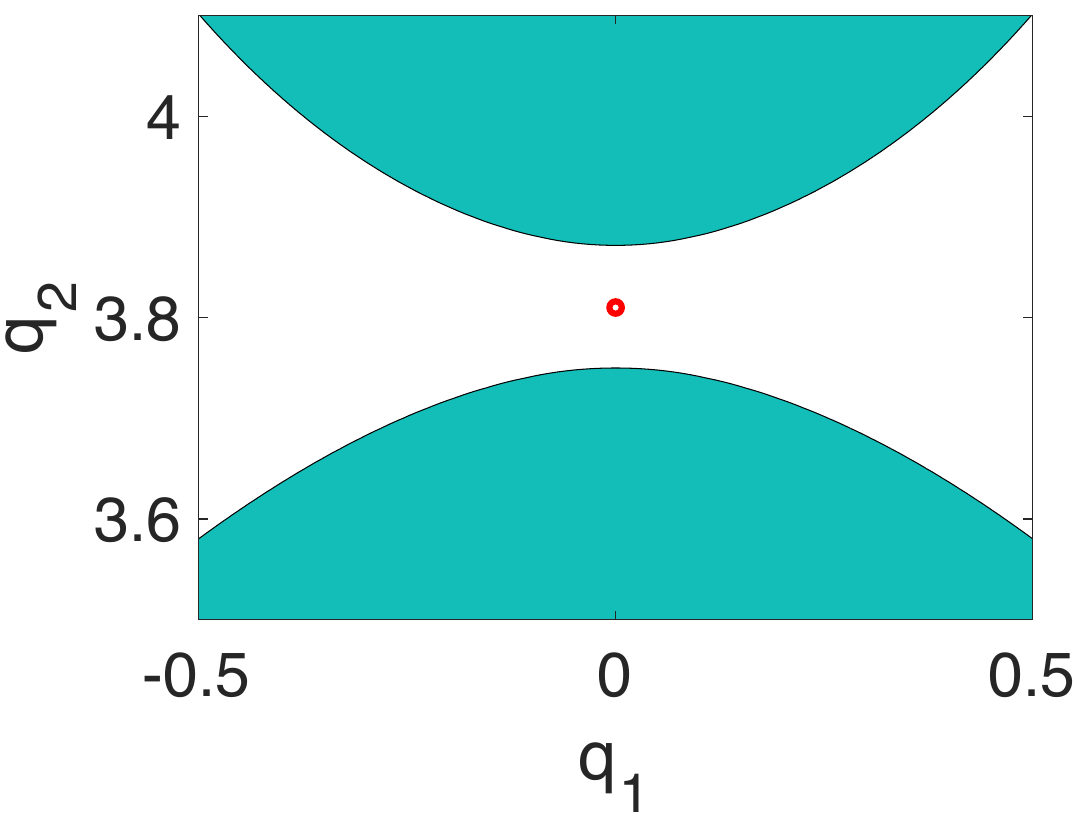} 
\includegraphics[width=.3\textwidth]{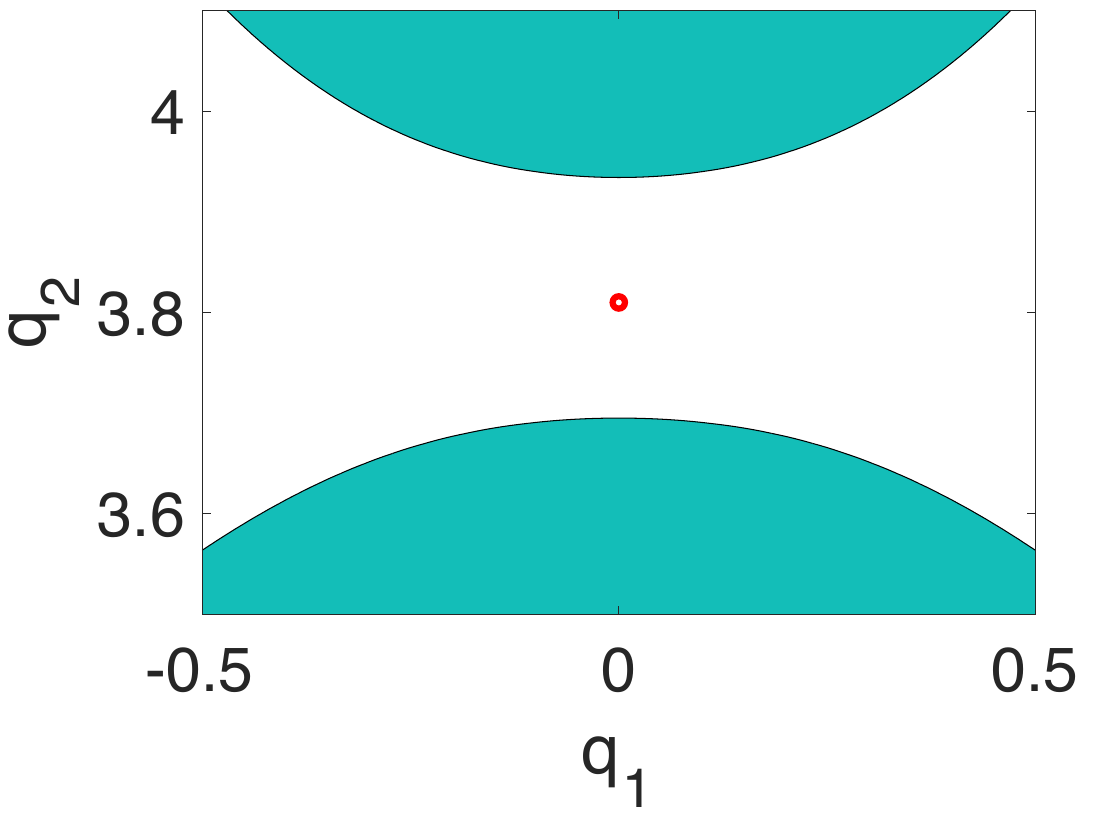}
\caption{\footnotesize The contours of potential energy $V(q_1,q_2)$ (black) separate the allowed (white) and forbidden (green) regions in the configuration space $(q_1,q_2)$ at total energy levels (left) $ E_1=E_{eq}$, (middle) $E_2>E_{eq}$, and 
 (right) $E_3>E_2$. The bottleneck around the equilibrium point (red circle) opens for $E>E_{eq}$}. 
\label{fig:bottleneck}
\end{figure}

In this section, we recall the main ideas involved in `tube dynamics' \cite{conley1968low,mcgehee1969some,koon2000dynamical}, and adapt the analysis to our system. Let $E_{eq}$ denote the energy of $saddle\times center$ equilibrium point $\mathcal{X}=\begin{pmatrix}0&0&q_2&0\end{pmatrix}^\intercal$. Since, the autonomous Hamiltonian system of Eqns. (\ref{eq:Ham_ode}) preserves energy, each trajectory lies on an isoenergetic 3D manifold in the 4D phase space. For energy levels $E$ slightly above $E_{eq}$, a bottleneck exists near the equilibrium in the $q_1-q_2$ plane, as shown in Fig. \ref{fig:bottleneck}. The trajectories can travel from left (i.e., $q_1<0$) to the right ($q_1>0$) region by passing through this bottleneck. The Jacobian matrix (Eq. \ref{eq:jac}) evaluated at $\mathcal{X}$ is of the form:
\begin{equation}
    A=\begin{pmatrix}
0 & -a & 0 & 0\\
-b & 0 & 0 & 0\\
0 & 0 & 0 & -c \\
0 & 0 & d & 0
\end{pmatrix},\label{eq:A}
\end{equation} where $a,b,c,d$ are all positive real numbers. As before, the linear dynamics of perturbation $\mathcal{Z}=\begin{pmatrix} q_{1l}&p_{1l}&q_{2l}&p_{2l}\end{pmatrix}^\intercal$ around $\mathcal{X}$ are given by $\dot{\mathcal{Z}}=A\mathcal{Z}$.  
The quadratic Hamiltonian corresponding to this linearized system is $H_l(q_{1l},p_{1l},q_{2l},p_{2l})=0.5(bq_{1l}^2-ap_{1l}^2-dq_{2l}^2-cp_{2l}^2).$ The  corresponding energy $E_l=-H_l$ is an invariant of the  linearized dynamics.
The eigenvalues of $A$ are $\lambda_{1,2}=\pm\sqrt{ab}$, and $\lambda_{3,4}=\pm i \sqrt{cd}$. The eigenvector matrix can be taken to be \begin{equation}\label{eq:T}
T=\begin{pmatrix}
    \sqrt{\dfrac{a}{a+b}}& \sqrt{\dfrac{a}{a+b}}&0&0\\
     -\sqrt{\dfrac{b}{a+b}} &\sqrt{\dfrac{b}{a+b}}&0&0\\
     0 &0& \sqrt{\dfrac{c}{c+d}} &0\\
     0&0&0&-\sqrt{\dfrac{d}{c+d}}
\end{pmatrix}.    
\end{equation}

The coordinates $\mathcal{Y}=\begin{pmatrix} \zeta&\eta&\rho_1&\rho_2\end{pmatrix}^\intercal$ in the eigenvector basis are given by the relation $\mathcal{Z}=T\mathcal{Y}$, and the corresponding linear dynamics $\dot{\mathcal{Y}}=T^{-1}AT\mathcal{Y}$ are:

\begin{equation}
 \begin{split}    
    &\dot{\zeta}=\lambda\zeta,\\&
    \dot{\eta}=-\lambda\eta,\\&
    \dot{\rho_1}=-\nu\rho_2,\\&
    \dot{\rho_2}=\nu\rho_1.
\label{eq:ham_lin_ebasis}
\end{split}
\end{equation}
The energy invariant in the new coordinates is $E_l(\zeta,\eta,\rho_1,\rho_2)=-a_1\zeta\eta+a_2(\rho_1^2+\rho_2^2)$, where $a_1=\dfrac{2ab}{a+b}$, and $a_2=\dfrac{0.5cd}{c+d}$. Since the dynamics of $(\zeta,\eta)$ and $(\rho_1,\rho_2)$ in Eqs. \ref{eq:ham_lin_ebasis} are decoupled, the system also possesses two additional invariants, namely $\zeta\eta$ and $\rho_1^2+\rho_2^2$.

Consider the region $\mathcal{R}$ defined by two constraints $E_l=\epsilon_1$ (fixed energy), and $|\zeta+\eta|\leq C$, where both $\epsilon_1$ and $C$ are positive. Rewriting the energy equation as 
\begin{equation}
    \dfrac{a_1}{4}(\zeta-\eta)^2+a_2(\rho_1^2+\rho_2^2)=\epsilon_1+\dfrac{a_1}{4}(\zeta+\eta)^2, \label{eq:ener_sphere}
\end{equation} we note that for $\zeta+\eta$ fixed, Eq. \ref{eq:ener_sphere} describes a topological 2-sphere (geometrically an ellipsoid) in the 4D phase space. Hence $\mathcal{R}$ has the topology $S^2\times I$. 

\begin{figure}[h!]
    \centering
    \includegraphics[width=.99\textwidth]{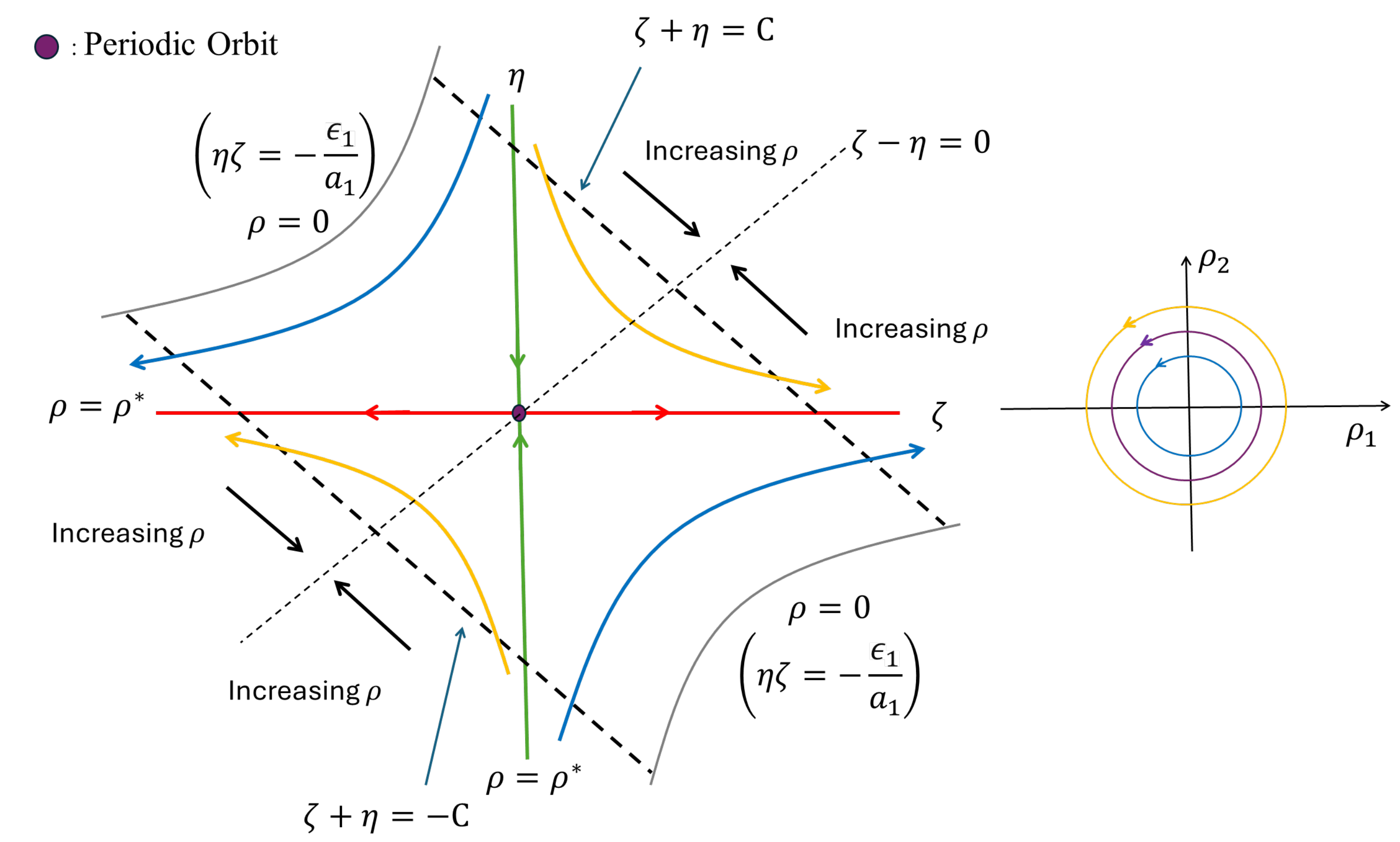}
    \caption{\footnotesize Phase space geometry in the $saddle\times center$ case. Projection of region $\mathcal{R}$ and the trajectories of the linearized Hamiltonian equations on the (left) $\zeta-\eta$ plane, and (right) $\rho_1-\rho_2$ plane. The periodic orbit (purple) projects to the origin on the former, and to a circle of radius $\rho^*$ on the latter. All other trajectories travel on cylinders, and their projections are hyperbolas and circles on the two planes, respectively. The trajectories on cylinders with radius $\rho<\rho^*$ (blue) transit between the bounding spheres, while those on cylinders with $\rho>\rho^*$ (yellow) return back to the originating bounding spheres. Trajectories on the red and green cylinders are asymptotic to the periodic orbit in negative and positive time, respectively. }
    \label{fig:R_etazeta}
\end{figure}

Fig. \ref{fig:R_etazeta} shows the projection of the $\mathcal{R}$ on the $\zeta-\eta$ and $\rho_1-\rho_2$ planes. Each point $(\zeta,\eta)$ in this projection corresponds to a circle in the $(\rho_1,\rho_2)$ plane, with radius $\rho=\sqrt{\rho_1^2+\rho_2^2}$ given by Eq. \ref{eq:ener_sphere}. There exists a periodic orbit (in the linear system), and it projects to the origin in the $\zeta-\eta$ plane. The $\zeta-\eta$ projection is bounded by the two dashed lines $\zeta+\eta=-C$ and $\zeta+\eta=C$, connecting the two boundary hyperbolas that correspond to $\zeta\eta=-\dfrac{\epsilon_1}{a_1}$ (and hence, $\rho=0$). The two 2-spheres in $\mathcal{R}$ corresponding to the boundary lines are called `bounding spheres' since they form the 2D boundary of the 3D region $\mathcal{R}$ in the 4D phase space. Each trajectory (except the periodic orbit) appears as a hyperbola ($\zeta\eta=\text{constant}$) in this projection, since $\zeta\eta$ is an invariant of motion as mentioned above. The following list relates the various types of trajectories in $\mathcal{R}$ and their $\zeta-\eta$ plane projections: \begin{enumerate}
    \item $\eta\zeta=0$ (axes): These are cylinders of radius $\rho^*=\sqrt{ \dfrac{\epsilon_1}{a_2}}$ in $\mathcal{R}$, corresponding to trajectories asymptotic to the periodic orbit in positive ($\eta=0$) or negative ($\zeta=0$) time. 
    \item $\eta\zeta<0$ (hyperbolas in second/fourth quadrants): These are trajectories lying on cylinders with $\rho<\rho^*$, and go from one bounding sphere to another.
    \item $\eta\zeta>0$ (hyperbolas in first/third quadrants): These are trajectories lying on cylinders with $\rho>\rho^*$. These cylinders begin and end at the same bounding sphere. 
\end{enumerate}

The above discussion is based on linear dynamics in the neighborhood of the equilibrium. The main takeaway is that only those trajectories that lie on cylinders with $\rho<\rho^*$, i.e. inside the tubes built up of orbits asymptotic to the periodic orbit (of the linear system), can transit from one boundary sphere to another. 

\begin{figure}[h!]
    \centering
   \includegraphics[width=.75\textwidth]{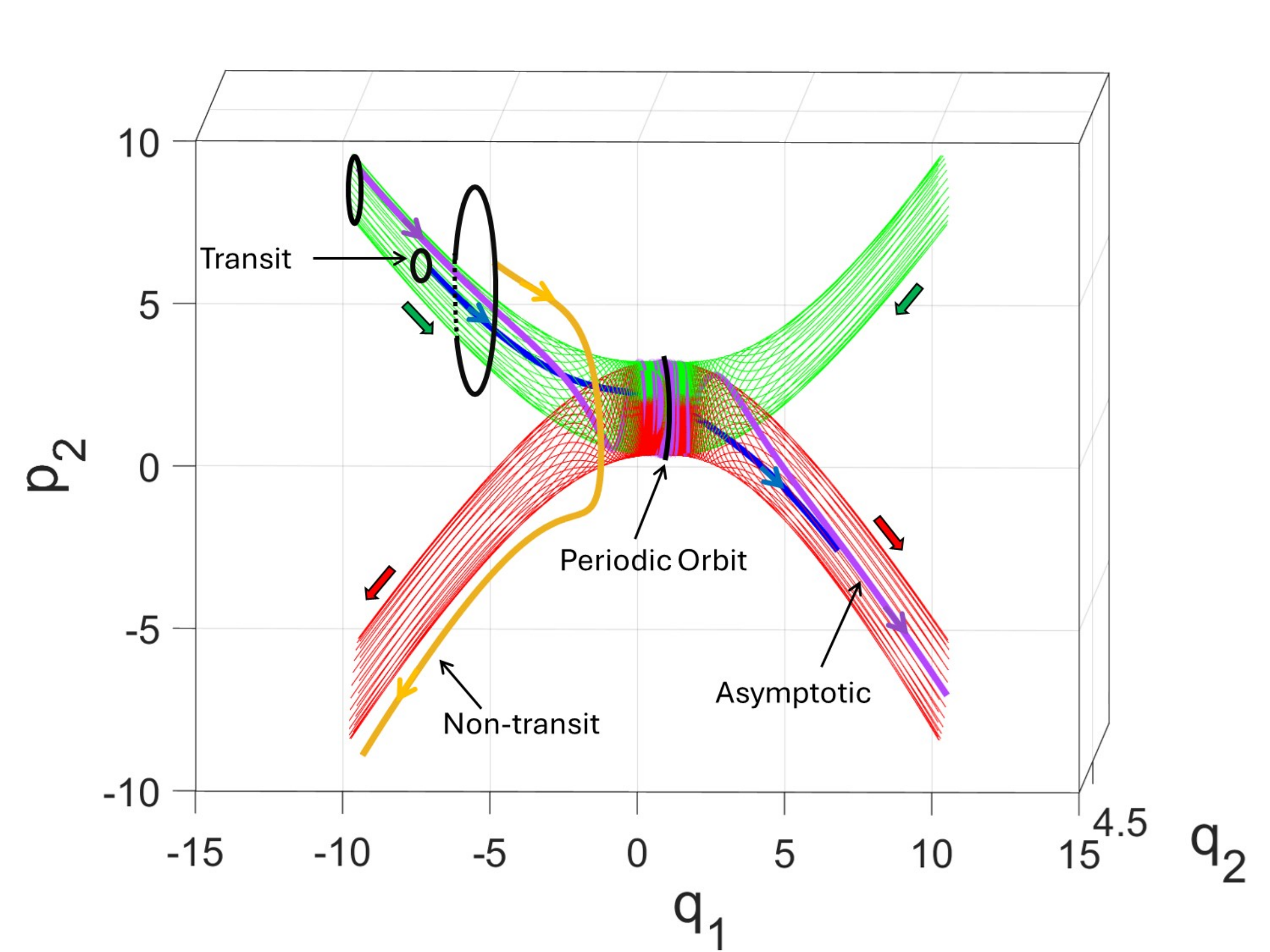}
    \caption{\footnotesize The 3D phase space geometry of the nonlinear Hamiltonian ODEs (restricted to a fixed energy level) near a $saddle\times center$ equilibrium point. The energy level $E$ is slightly above that of the equilibrium. The tube-shaped stable (green) and unstable (red) manifolds of the periodic orbit form  barriers to transport in this system. Analogous to the linear picture of Fig. \ref{fig:R_etazeta}, the trajectories starting inside (blue) the stable tube transit across, while those starting outside (yellow) the tube do not. The transiting trajectories that stay (approximately) on the tubes (purple) are referred to as `asymptotic' in the main text. }
    \label{fig:tubes}
\end{figure}

For sufficiently small positive values of excess energy $\epsilon_1=E-E_{eq}$, a family of periodic orbits of the nonlinear Hamiltonian ODE (Eqs. \ref{eq:Ham_ode}) is guaranteed to exist around $\mathcal{X}$ \cite{koon2000dynamical}. In this case, the above described qualitative picture of phase space near the equilibrium persists in the nonlinear system. For the periodic orbit of the nonlinear system $\mathcal{P}(\tau):[0,1]\rightarrow \mathbb{R}^4$ (satisfying $\mathcal{P}(0)=\mathcal{P}(1)$), there exist 2D stable ($W^s_{p.o.}$) and unstable ($W^u_{p.o.}$) manifolds:
\begin{align}
    W^s_{p.o.}=\{Z\in\mathbb{R}^4|\lim_{t\rightarrow\infty}\phi_t(Z)\rightarrow \mathcal{P}\},\\
    W^u_{p.o.}=\{Z\in\mathcal{R}^4|\lim_{t\rightarrow -\infty}\phi_t(Z)\rightarrow \mathcal{P}\}.
\end{align}
These manifolds are locally diffeomorphic to cylinders. Only those trajectories that are inside these tube-like stable/unstable manifolds can transit, see Fig. \ref{fig:tubes}. This is a global result since these 2D invariant manifolds are codimension-1 in the 3D phase space (of fixed energy), and hence trajectories cannot cross the surface of these tubes. For details on computation of periodic 
orbits and their invariant manifolds, we refer the reader to \cite{koon2000dynamical}.

\section{Bifurcations in the Hamiltonian BVP : Phase space analysis and numerical continuation}

In this section, we discuss the numerical solutions of the Hamiltonian BVP (Eqs. \ref{eq:Ham_ode}, \ref{eq:ic_fc_ham}) as the time horizon $T$ is varied, and interpret the solution structure using the geometry of the phase space for the two cases discussed in the previous section. In both cases, we use the following initial and final conditions:
\begin{equation}
 \begin{split} 
    &q_1(0)=-10, q_2(0)=4.5,\\&
    q_1(T)=10, q_2(T)=4.5.
     \end{split} 
     \label{eq:ic_bc}
\end{equation}

The solutions are computed using a combination of the boundary value problem solver BVP4c \cite{shampine2000solving}, and the Computational Continuation Core (COCO) \cite{dankowicz2013recipes} in Matlab. 

\begin{figure}[hbt!]
\centering
\begin{tabular}{c c}
\includegraphics[width=.5\textwidth]{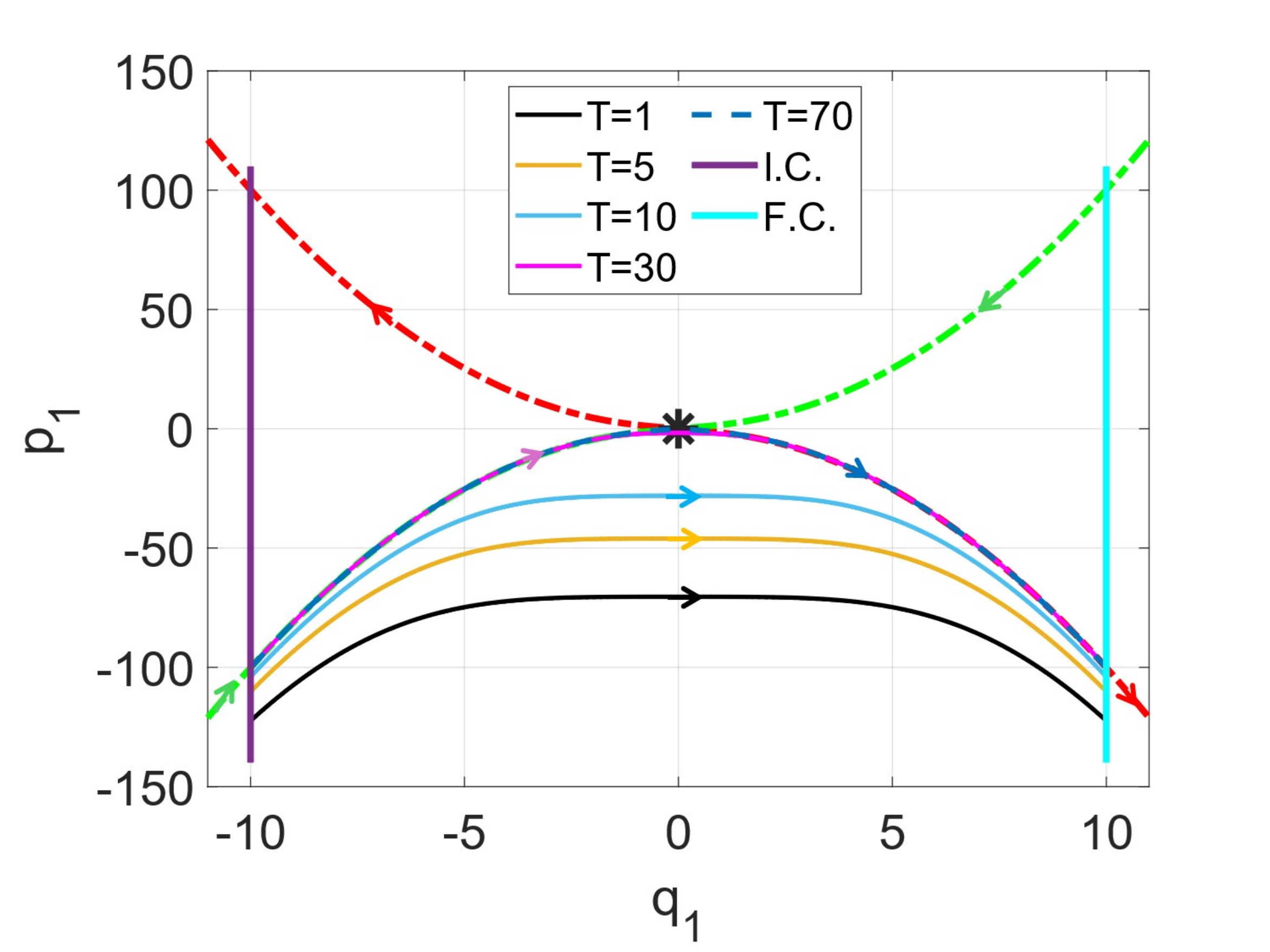} &\includegraphics[width=.5\textwidth]{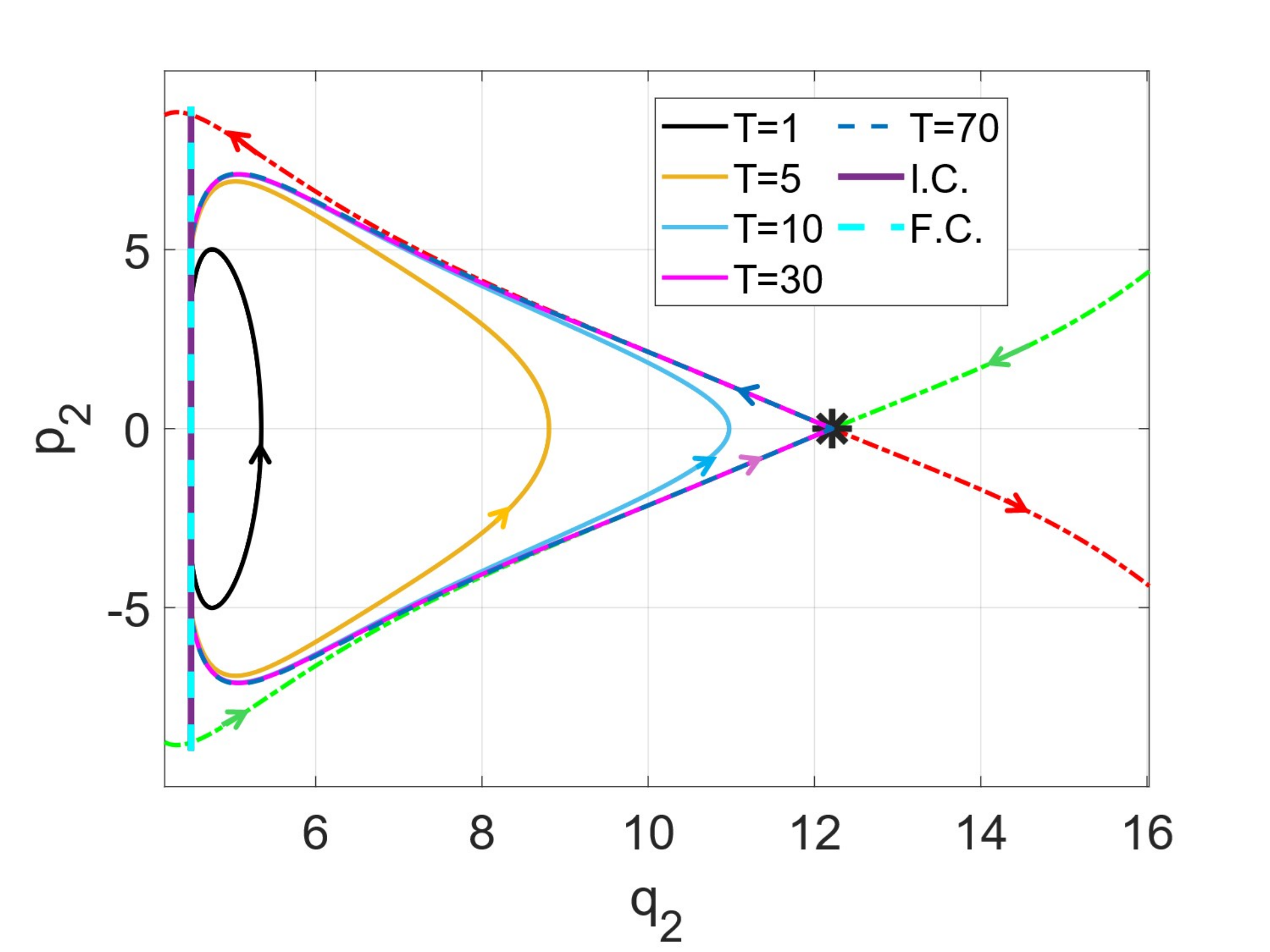} 
\\ (a) & (b)
\end{tabular}
\caption{\footnotesize Various solutions of the Hamiltonian BVP in the $saddle\times saddle$ equilibrium case shown in the (a) $q_1-p_1$, and (b) $q_2-p_2$ phase planes, along with the projections of the initial and final conditions. Also shown are the projections of the stable (green) and unstable (red) manifolds of the equilibrium. As $T$ is increased, the trajectories get closer to the invariant manifolds. 
}
\label{fig:bvp_ss_phase}
\end{figure}

\begin{figure}
    \centering
 \includegraphics[width=.5\textwidth]{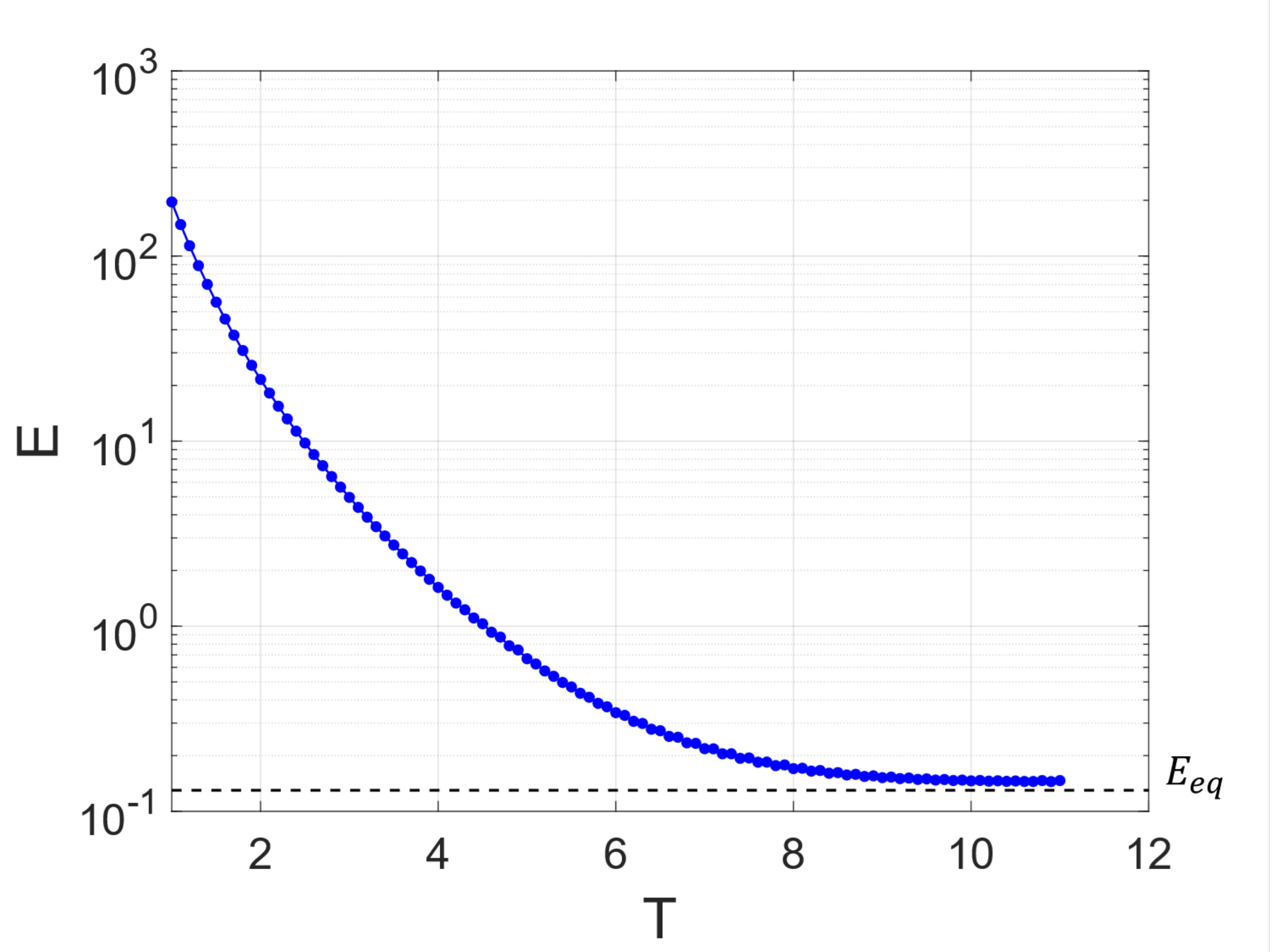}
    \caption{\footnotesize The energy (E) v/s time-horizon (T) diagram for the $saddle\times saddle$ case. As $T\rightarrow\infty$, the energy approaches that of the equilibrium.}
    \label{fig:bif_ss}
\end{figure}
\subsection{Hamiltonian BVP with a $saddle\times saddle$ type ergodic equilibrium}
We pick parameters values $\sigma=1,\mu=2,g=4,h=0,\alpha=1$, and $\epsilon=0.05$, in which case the system has a unique equilibrium $\mathcal{X}=(q_1 \;p_1\; q_2\; p_2)^{\intercal}=(0\; 0\; 12.21 \;0)^{\intercal}$. $\mathcal{X}$ is a $saddle\times saddle$ type equilibrium, with the eigenvector matrix of the form $T$ given by Eq. \ref{eq:T_ss}, where $a=0.5,b=1.1186,c=200$ and $d=229.74$. 

In Fig. \ref{fig:bvp_ss_phase}, several solutions of the BVP are shown via projections onto the $q_1-p_1$ and $q_2-p_2$ planes, along with the stable and unstable manifolds of the equilibrium. There is a unique trajectory for each value of $T$, determined by the geometry of the stable and unstable manifolds. With increasing $T$, the trajectories get progressively closer to these invariant manifolds. In the limit $T\to\infty$, the solutions converge to the invariant manifolds, while the energy $E$ approaches that of the equilibrium $E_{eq}$, as shown in Fig. \ref{fig:bif_ss}. The limiting behavior is exactly the ergodic regime of the MFG solution. Overall, the solutions of the coupled BVP with a $saddle\times saddle$ fixed point behave similarly to those of the uncoupled case considered in \cite{ullmo2019quadratic}.

\subsection{Hamiltonian BVP with a $saddle\times center$ type ergodic equilibrium}
\begin{figure}[h!]
    \centering
    \includegraphics[width=.995\textwidth]{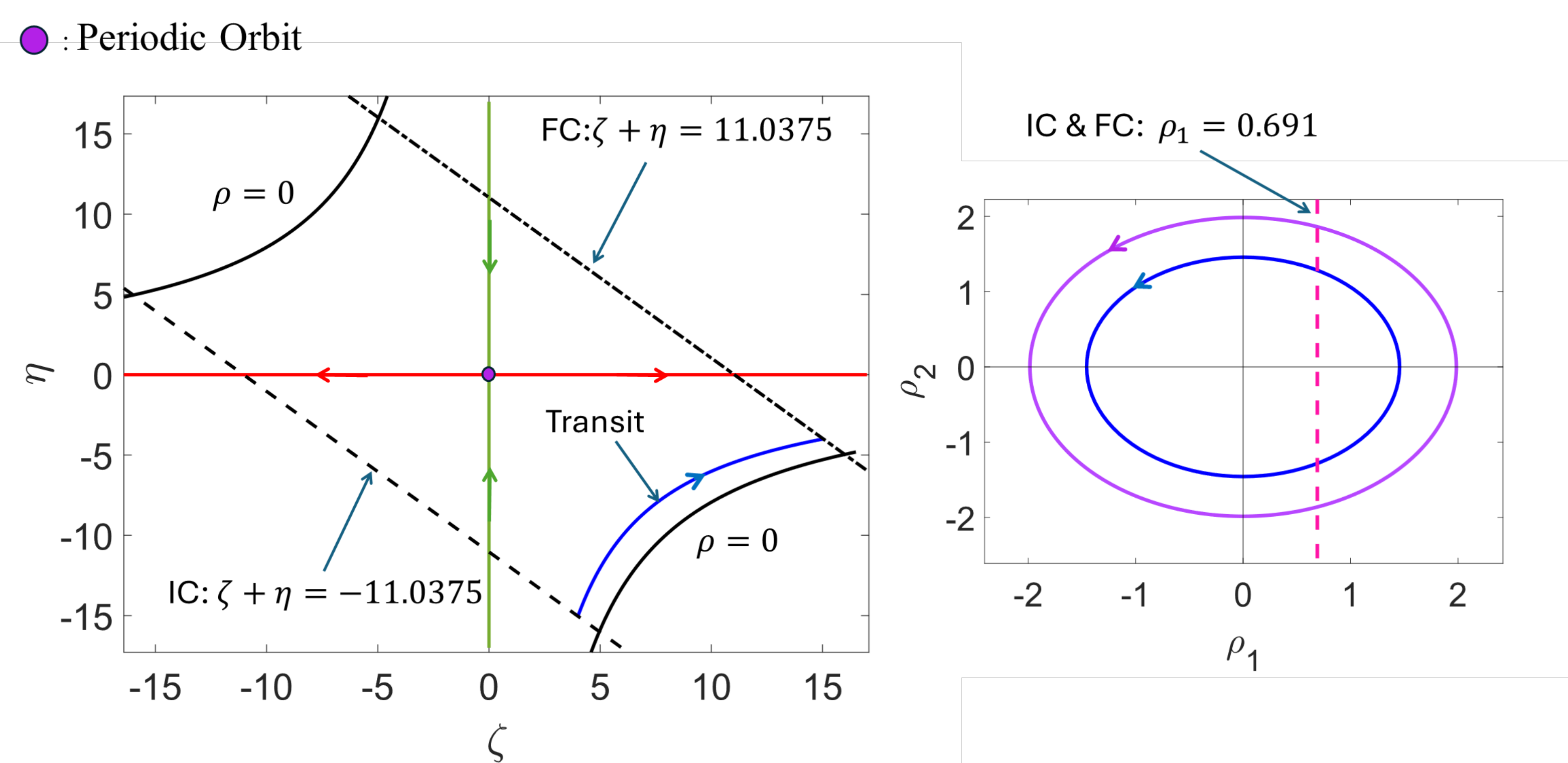}
    \caption{\footnotesize Projection of region $\mathcal{R}$ and a solution of the linearized BVP equations on the (left) $\zeta-\eta$ plane, and (right) $\rho_1-\rho_2$ plane. The BVP solutions (blue) must transit from the line of initial conditions (dash) to the line of final conditions (dash-dot) in the $\zeta-\eta$ plane, and hence lie on cylinders with $\rho<\rho^*$. }
    \label{fig:bvp_R}
\end{figure}

In this case, we pick the parameters $\sigma=1, \mu=2, g=4, \alpha=3, h=0$ and $\epsilon=0.05$, which results in a $saddle\times center$ type equilbrium point $\mathcal{X}=(q_1 \;p_1\; q_2\; p_2)^{\intercal}=(0\; 0\; 3.81 \;0)^{\intercal}$. The Jacobian matrix $A$ evaluated at $\mathcal{X}$ is of the form given in Eq. \ref{eq:A}, with $a=0.5 ,b=0.109,c=200$ and $d=0.946$, and eigenvalues $\pm 0.233$ and $\pm 13.8i$.  The system possesses another equilibrium point (with $q_2\approx 100$) that is irrelevant for the chosen initial and final boundary conditions, so we focus our discussion on the region around $\mathcal{X}$.

Similar to the discussion in Sec. \ref{subsec:phase_sc}, we first analyze the linearized system by restricting to an isoenergetic 3D region $\mathcal{R}$ with $E_l=\epsilon_1$. Using the transformation matrix $T$ (Eq. \ref{eq:T}) and the relation $\mathcal{Z}=T\mathcal{Y}$, we obtain $q_1=0.906(\zeta+\eta)$. Hence, the set of initial and final conditions defined by $q_1(0)=-10$ and $q_1(T)=10$ projects on the $\zeta-\eta$ plane parallel to the boundary lines, see Fig. \ref{fig:bvp_R} (left). Similarly, the other set of initial and final conditions defined by $q_2(0)=q_2(T)=4.5$ is shown in the $\rho_1-\rho_2$ plane in Fig. \ref{fig:bvp_R}(right).
\begin{figure}[h!]
\centering
\includegraphics[width=.75\textwidth]{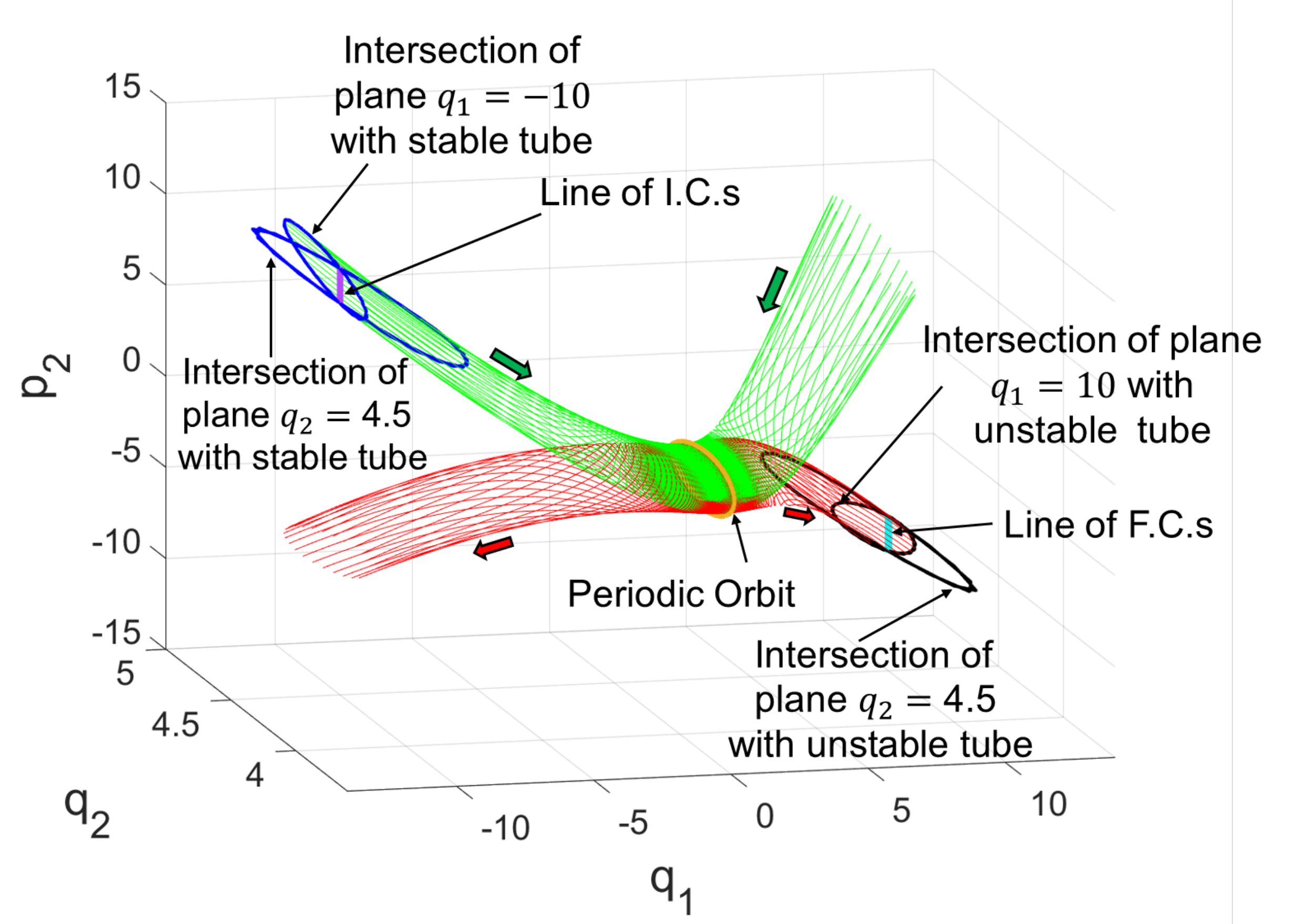} 
\caption{\footnotesize The 3D phase space geometry of the nonlinear Hamiltonian BVP (restricted to a fixed energy level) near a $saddle\times center$ equilibrium point, showing the invariant manifolds (tubes) of the periodic orbit, and the lines of initial/final conditions. Each BVP solution starts at the line of initial conditions, travels inside/on the stable tube (green), and then switches to the unstable tube (red), ending on the line of final conditions.} \label{fig:bvp_tubes}
\end{figure}

As discussed in the Sec. \ref{subsec:phase_sc}, the only points that transit from the initial to the final condition line in the $\zeta-\eta$ plane are those with $\rho<\rho^*$. 
Simultaneously, the trajectory has to begin and end at the line of initial and final conditions in the $\rho_1-\rho_2$ plane. 

The corresponding picture in the 3D phase space of the nonlinear system is shown in Fig. \ref{fig:bvp_tubes}, where we fix energy slightly above that of the equilibrium. The two planes of initial conditions $q_1=-10$, and $q_2=4.5$ intersect the stable (solid) tube in topological discs. The intersection of these two discs yields a line segment of feasible initial conditions for each level of energy. Similarly, the final conditions yield the feasible line segment upon intersection with the unstable tube. Each BVP solution with the prescribed energy level must begin on the starting line segment, and end on the final line segment in this 3D phase space. Each solution can be divided into three phases: 1). The arrival phase $0\leq t\leq t_a$, during which it travels from the initial condition to the neighborhood of the equilibrium, 2). The ergodic phase $t_a\leq t\leq t_a+\tau_{erg}$ during which the trajectory stays close to the equilibrium , and 3). The departure phase $t_a+\tau_{erg}\leq t\leq t_a+\tau_{erg}+t_d\approx T$, corresponding to the travel from the neighborhood of the equilibrium to the final condition. 

\begin{figure}[hbt!]
\centering
\includegraphics[width=.99\textwidth]{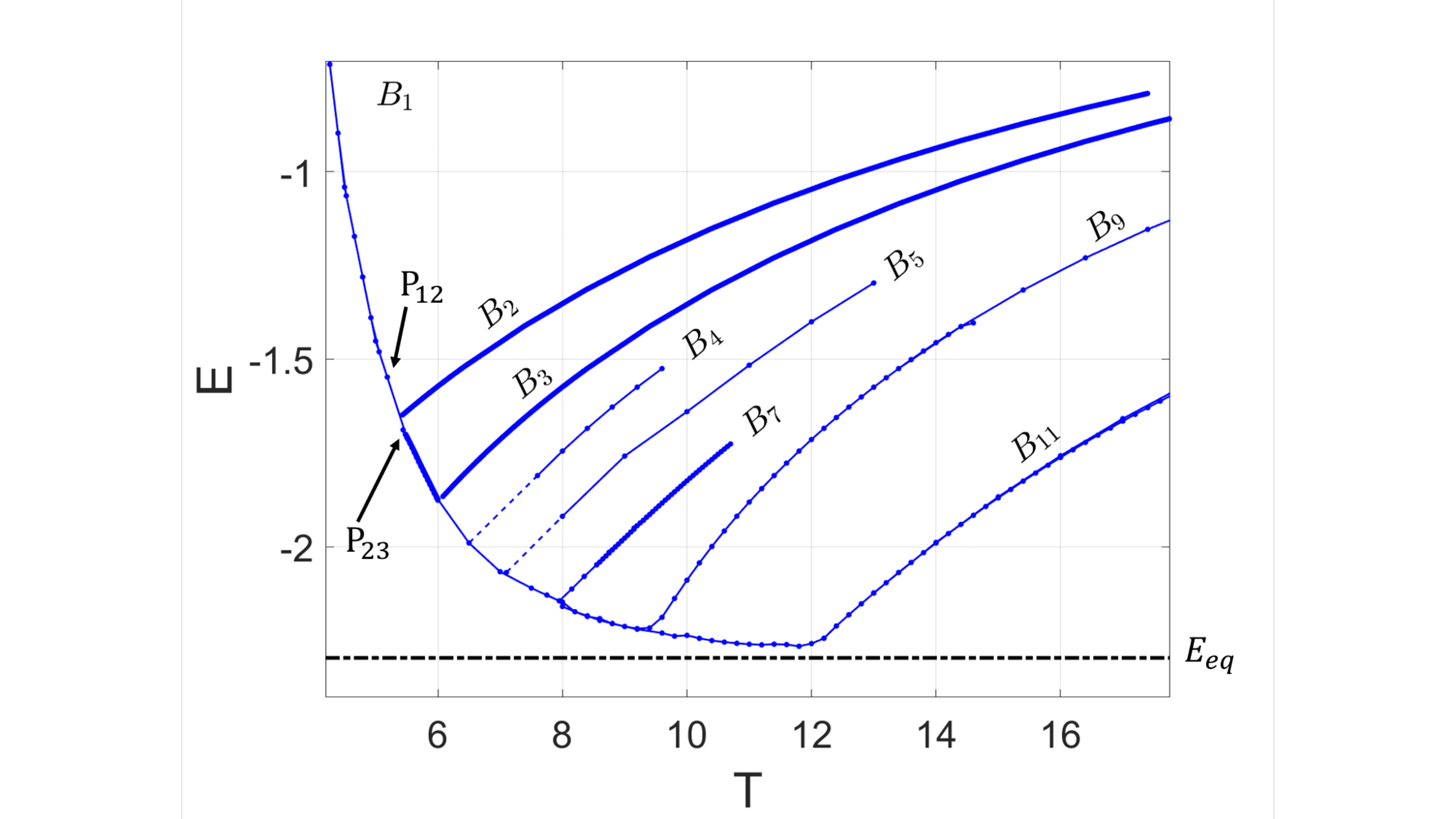}
\caption{\footnotesize A partial energy (E) v/s time-horizon (T) diagram for the $saddle\times center$ case, where solutions on branch $B_n$ have $n$ intersections with the $p_2=0$ plane, and undergo $n-1$ half-rotations during the ergodic phase. The $B_1$ branch is similar to the solution branch of the $saddle\times saddle$ case. The trajectories on branches $B_2$ through $B_{11}$ consist of asymptotic trajectories travelling (nearly) on the tubes. The trajectories lying on $B_1$ as well those on the segments connecting different branches, travel inside the tubes.}
\label{fig:bif_main}
\end{figure}

The bifurcation diagram for this system is shown in Fig. \ref{fig:bif_main}. The diagram contains several branches, and the system exhibits multiplicity of solutions (at fixed $T$) for $T>5.5$. Along the solution branches $B_2$ through $B_{11}$, the energy increases monotonically with $T$. Each solution on these branches is of asymptotic type, i.e. it travels (approximately) `on' the tubes rather that inside them. Hence, the time period of one rotation around the tubes along the trajectory (in the ergodic phase) is approximately equal to that of the periodic orbit ($t_p$), which in turn is fixed by the value of $E$. The $n$th bifurcation point, from which the $B_{n+1}$ branch originates, corresponds to the lowest energy level yielding $n$ half-rotations during the ergodic phase of the trajectory, i.e., $\tau_{erg}=0.5nt_p$. Along $B_{n+1}$, the time period $t_p$ of the periodic orbit as well as the tube radius increase with $E$. Most of the concomitant increase in $T$ is due to increase in $\tau_{erg}$, and it is such that the rotation remains fixed, i.e., $\dfrac{\tau_{erg}}{t_p}=0.5n$, for all points on that branch. Hence, all solutions on a single branch $B_n (2\leq n\leq 11)$ have the same topology. For a fixed time-horizon $T>5.5$, a solution can travel on a bigger tube (higher $E$) and do fewer rotations, or on smaller tube(s) (lower $E$) with more rotations.

\begin{figure}
\begin{tabular}{c c}\includegraphics[width=.55\textwidth]{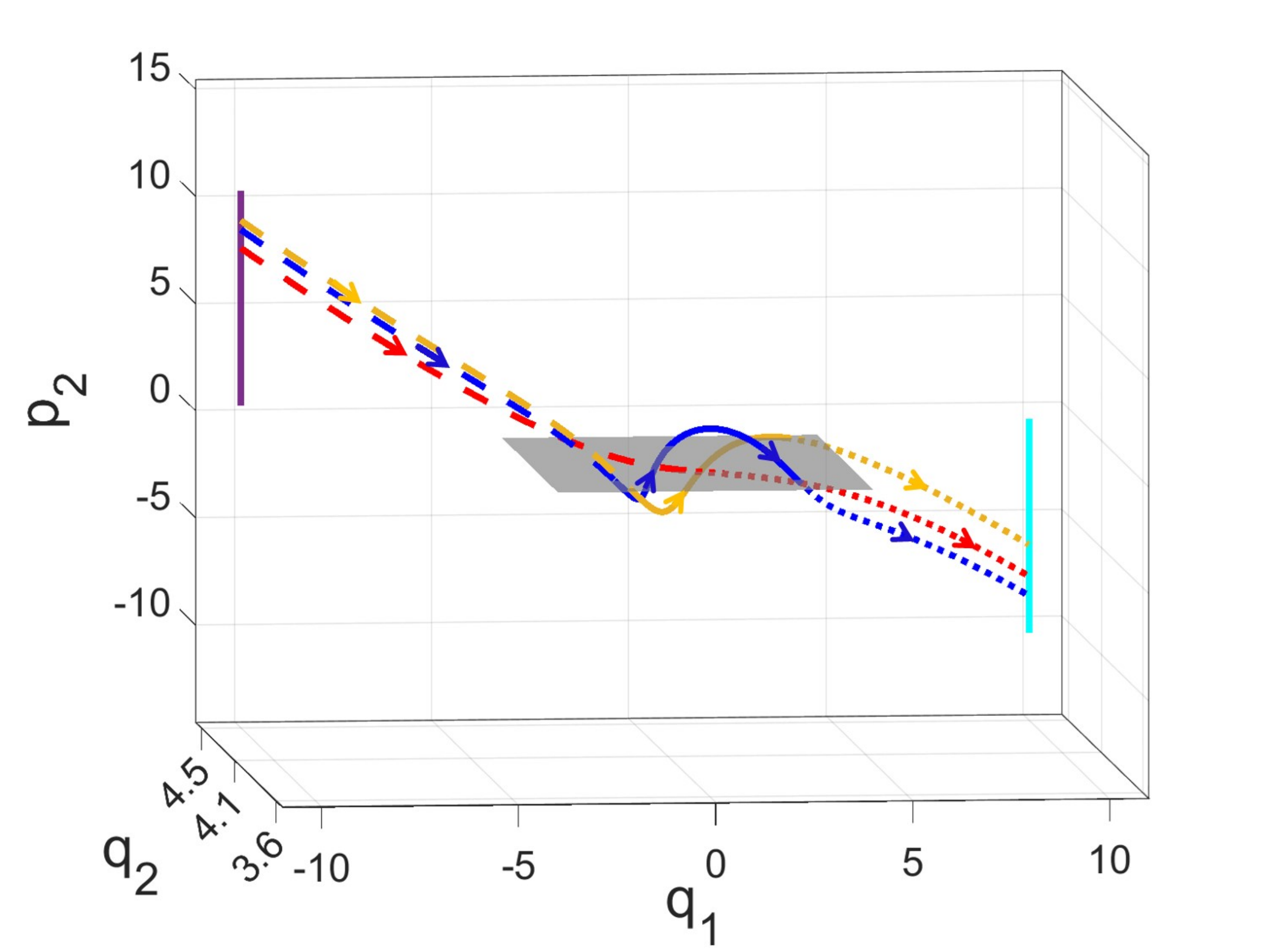} &\hspace{-.3in}\includegraphics[width=.55\textwidth]{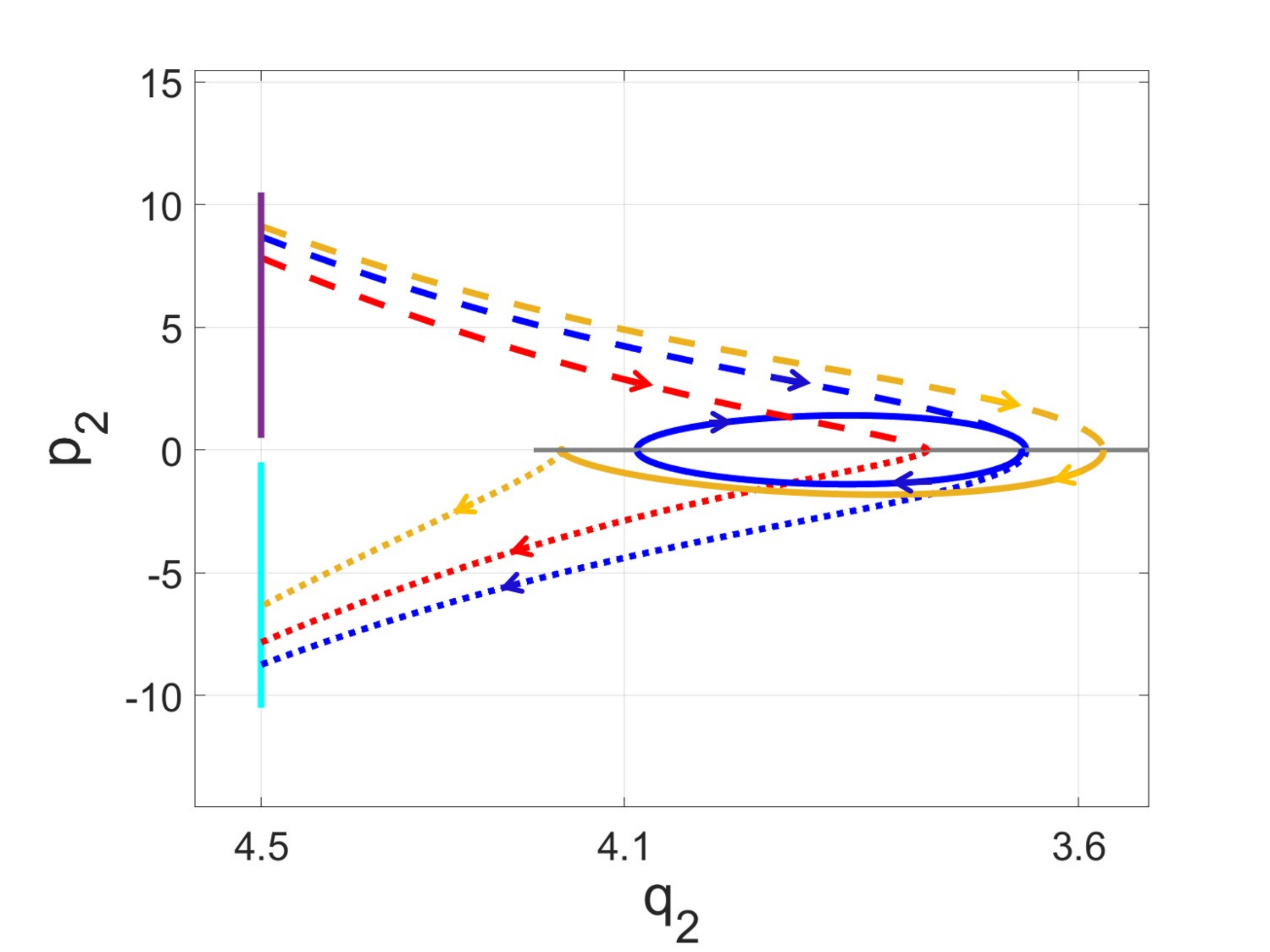} 
\\ (a) & (b)
\end{tabular}
\caption{\footnotesize Trajectories from $B_1$ (red), $B_2$(yellow) and $B_3$(blue) branches shown in the phase space, along with the lines of initial (purple) and final (cyan) conditions, as well as the plane $p_2=0$. The arrival (dash), and departure (dot) phases of all three trajectories are similar. In the ergodic phase (bold) the $B_1$ trajectory has no rotations, the $B_1$ trajectory completes half a rotation while the $B_2$ trajectory completes a full rotation.}
\label{fig:B1B2B3}
\end{figure}

\begin{figure}[h!]
\centering
\begin{tabular}{c c}
\hspace{-0.2in}
\includegraphics[width=.5\textwidth]{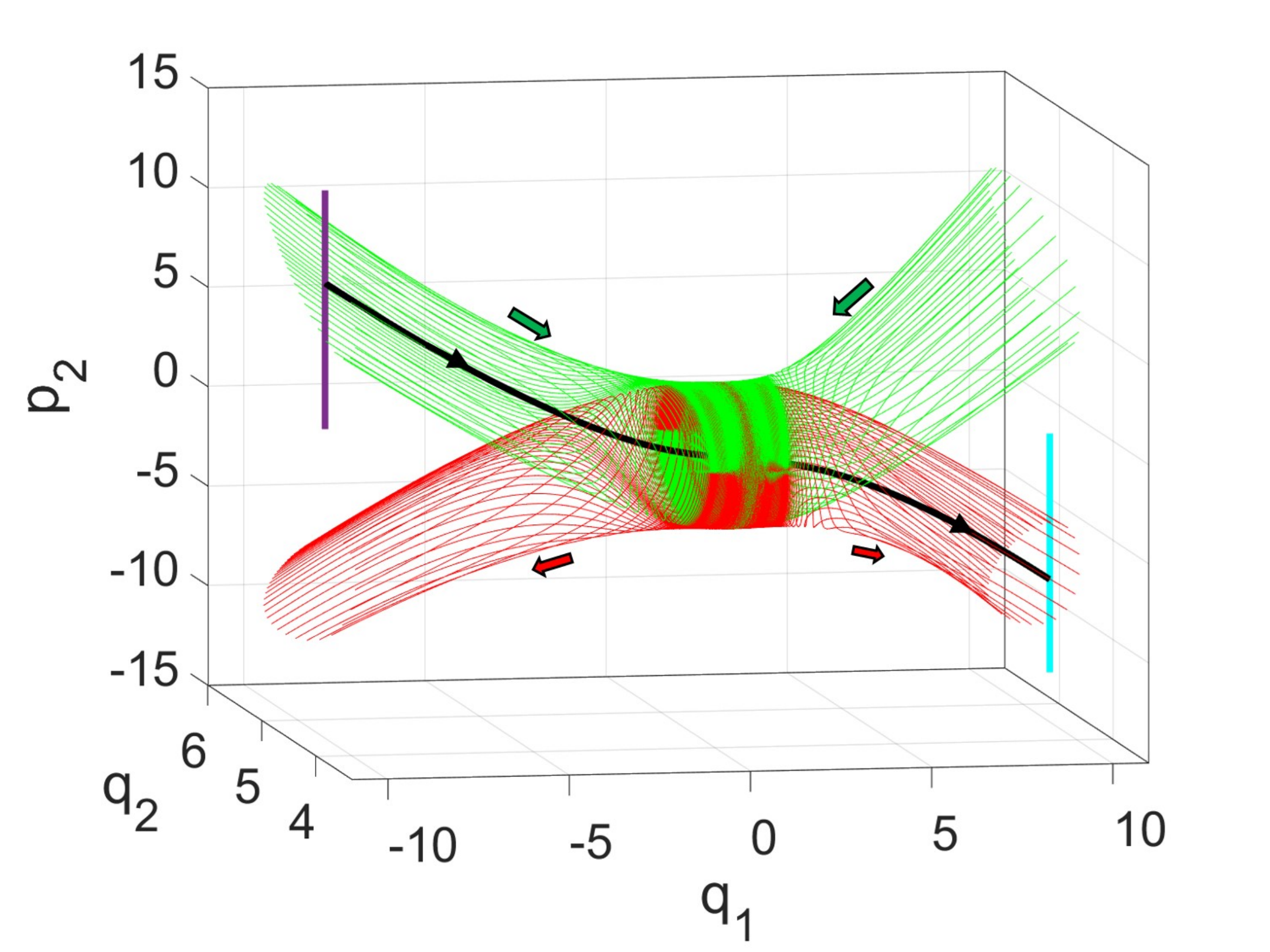} & 
    \includegraphics[width=.5\textwidth]{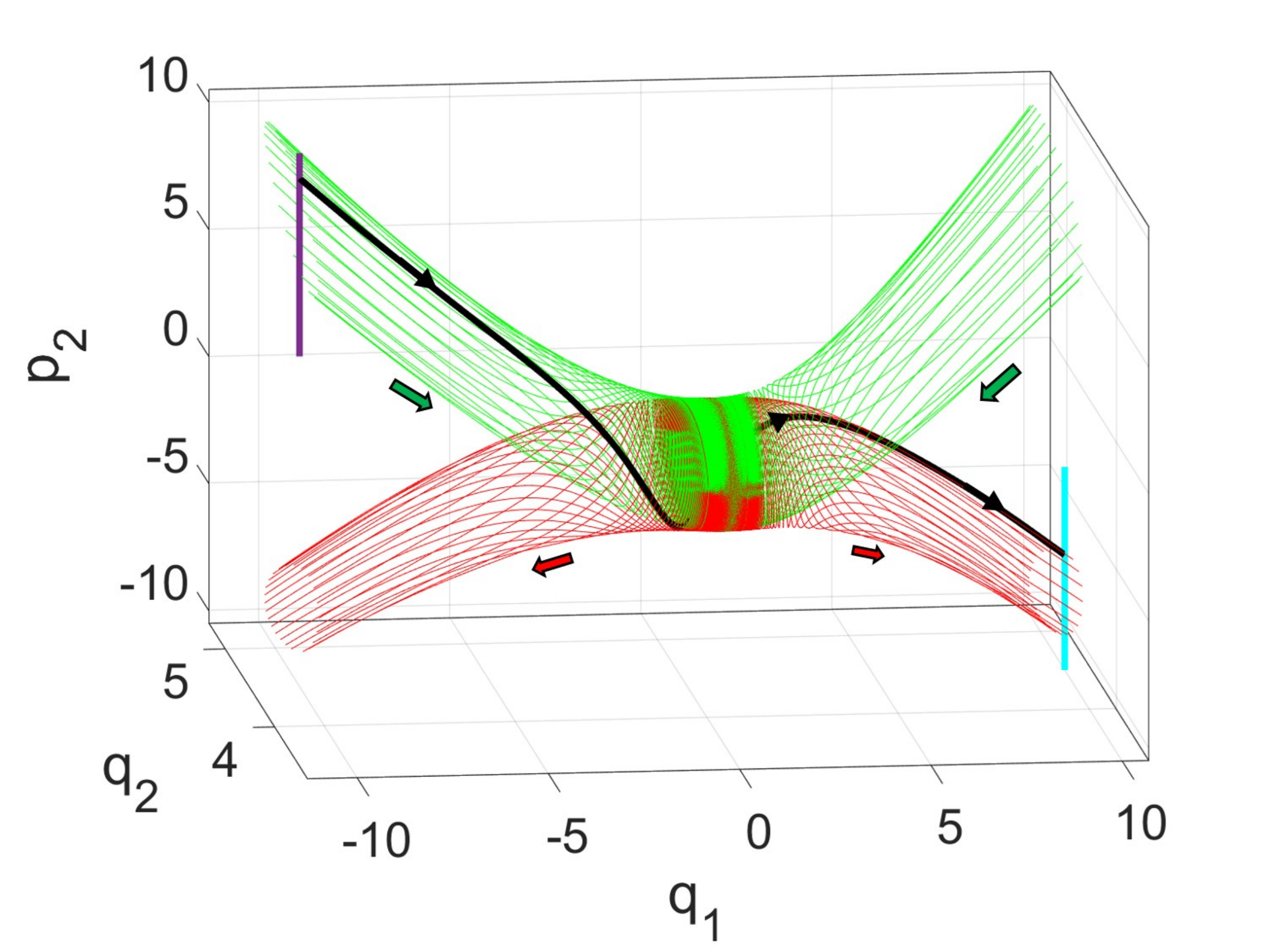} \\ (a)& (b)\\ 
\includegraphics[width=.5\textwidth]{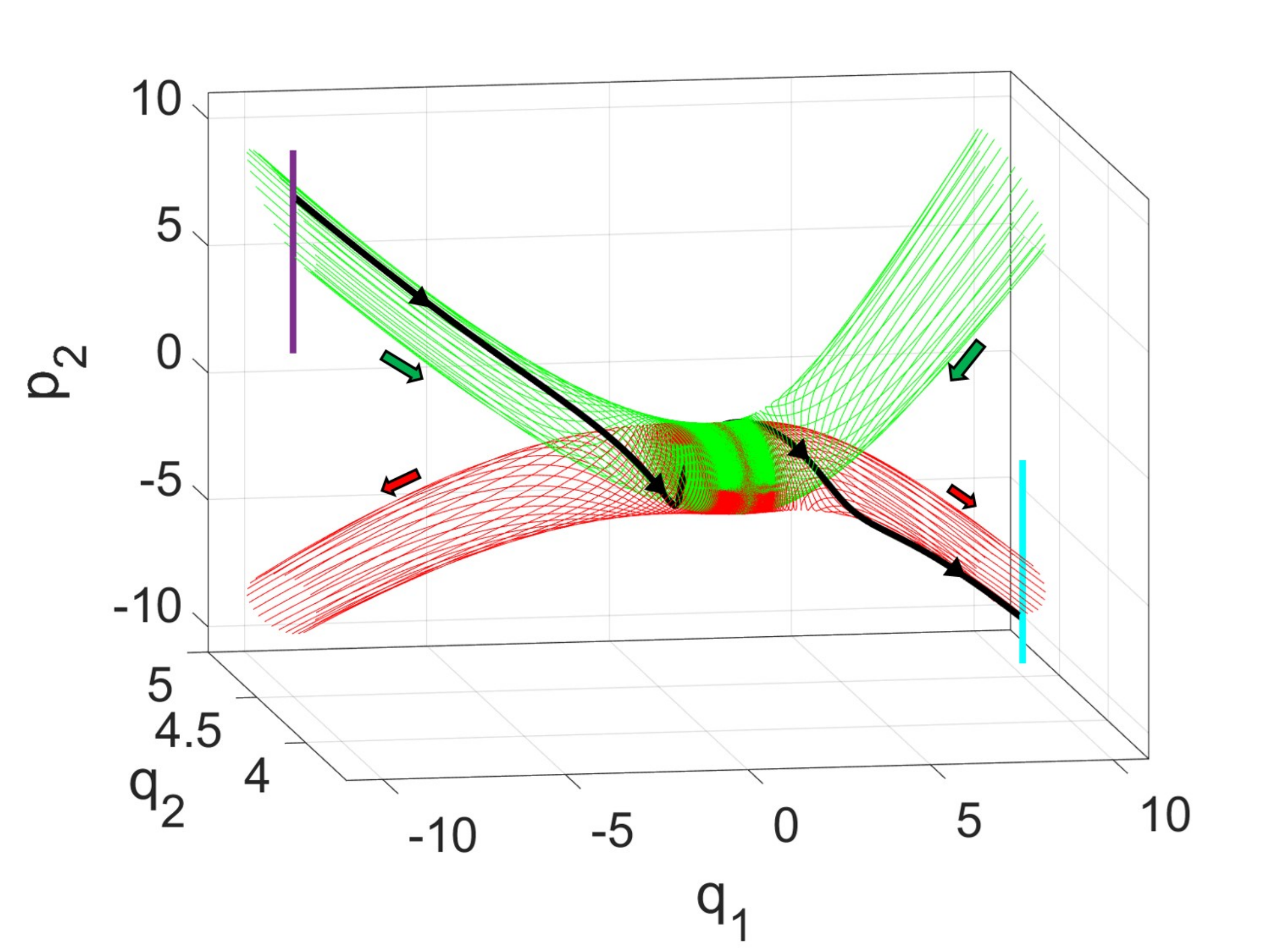} & 
\includegraphics[width=.5\textwidth]{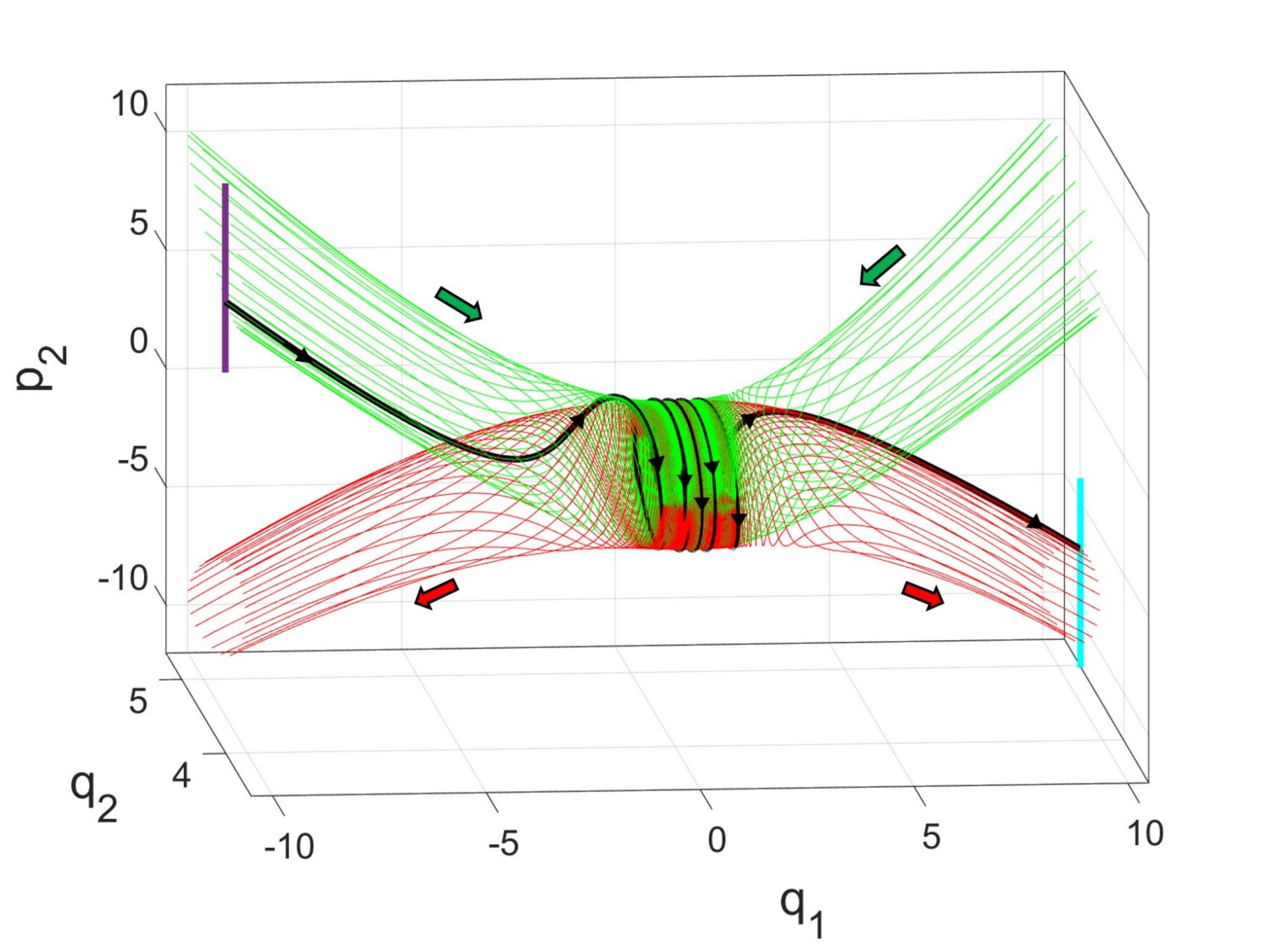} 
\\ (c)& (d)\\ 
\end{tabular}
\caption{\footnotesize The three trajectories of Fig. \ref{fig:B1B2B3} along with the stable and unstable tubes at the energy level of each trajectory. The (a) $B_1$, (b) $B_2$, and (c) $B_3$ trajectories complete zero, one-half and one rotations in the ergodic phase, respectively. (d) For comparison, a trajectory from the $B_{11}$ branch, completing five rotations in the ergodic phase. All trajectories except the one on the $B_1$ branch travel on the tubes.} \label{fig:B1B2B3_tubes}
\end{figure}

The topological origin of different branches is further evident in Fig. \ref{fig:B1B2B3}. This figure shows three trajectories from the $B_1, B_2$ and $B_3$ branches along with the lines of initial and final conditions. The arrival and departure phase of all trajectories correspond approximately to the segments connecting the initial condition line to the $p_2=0$ plane, and the segments from that plane to the line of final conditions, respectively. The trajectory on $B_1$ has only 1 intersection with the $p_2=0$ plane in the ergodic phase, while that on $B_2$ completes half a rotation, resulting in two intersections, and the one on $B_3$ completes a full rotation resulting in three intersections. The same three trajectories are shown in the 3D phase space in Fig. \ref{fig:B1B2B3_tubes}. Our attempts at numerical continuation of the branch $B_2$ failed at an energy level very close to the value above which a periodic orbit could not be found in the nonlinear system. Although we did not continue upto failure when computing branches $B_3$ through $B_{11}$, we expect a similar conclusion to hold in those cases too.

\begin{figure}[hbt!]
\centering
\includegraphics[width=.99\textwidth]{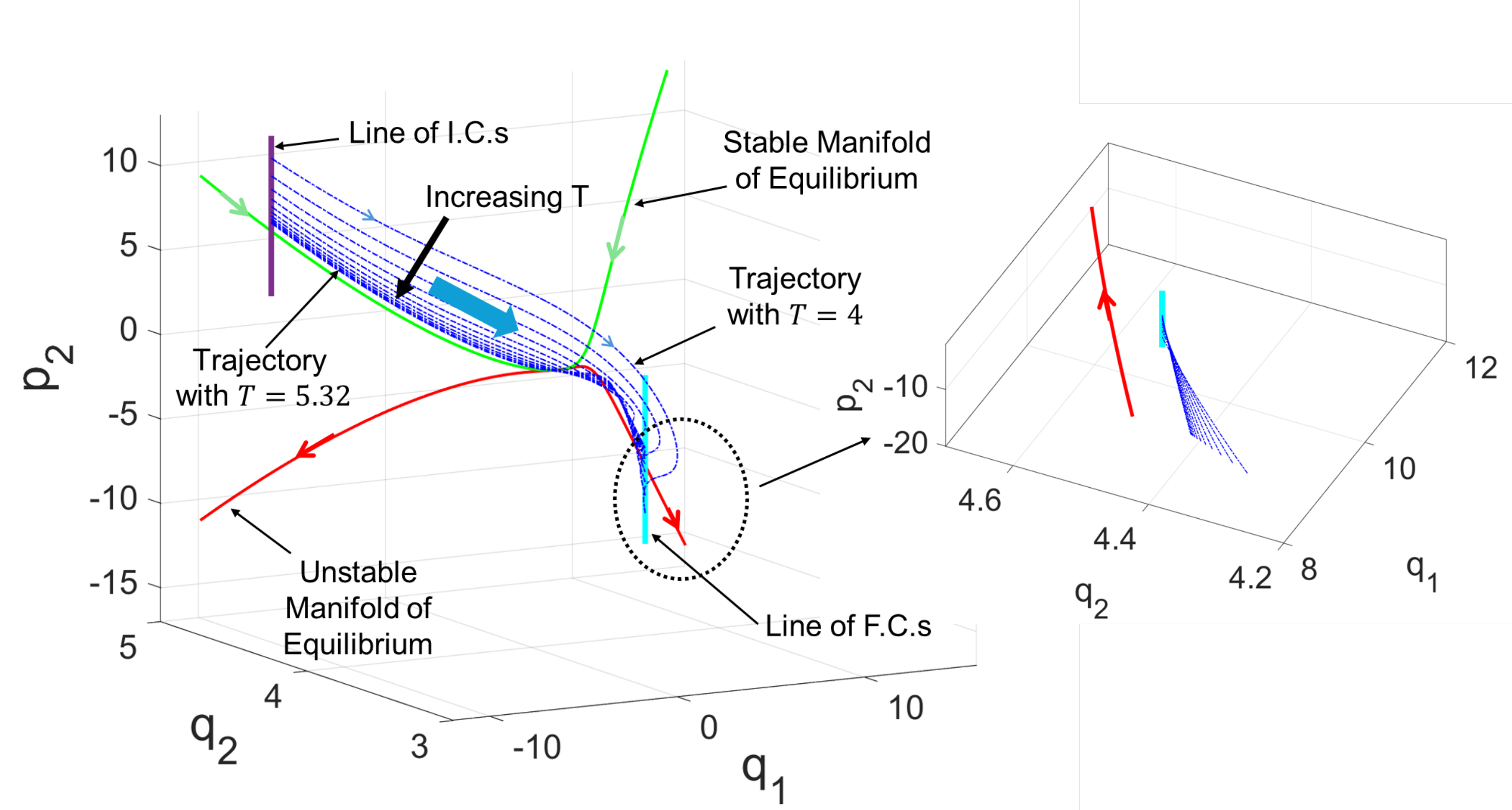}
\caption{\footnotesize 3D phase portrait showing several BVP solutions on the $B_1$ branch, along with the stable (green) and unstable (red) manifolds of the equilibrium. As the time-horizon $T$ is increased, the distance between trajectories and the manifolds asymptotes down to a finite value (at $T\approx 5.32$, labeled $P_{12}$ in Fig. \ref{fig:bif_main}). (Inset) The unstable manifold does not intersect the line of final conditions, hence the the radius of cylinders of the BVP trajectories cannot become arbitrarily small.}
\label{fig:B1_3D}
\end{figure}

The $B_1$ branch is similar to the only solution branch found in the $saddle\times saddle$ case in the previous section.  As $T$ is increased, trajectories get progressively closer to the stable and unstable manifolds of the equilibrium point, see Fig. \ref{fig:B1_3D}. However, these 1D invariant manifolds of the equilibrium point do not intersect the line of final condition, as shown in the inset of Fig. \ref{fig:B1_3D}. This implies a lower bound on the radius of cylinders on which the BVP trajectories can lie in order to satisfy the final condition. The trajectory corresponding to $T \approx 5.32$ on that branch hits this lower limit on the cylinder size, hence the branch terminates at this point, labelled $P_{12}$ on Fig. \ref{fig:bif_main}. Although we couldn't converge onto solutions between $P_{12}$ and the lowest energy point on $B_2$, we conjecture that they behave analogous to those on the segment starting at $P_{23}$ and ending on lowest energy solution on $B_3$ branch. 
These latter trajectories have same topology as $B_3$ but travel inside the tubes (rather than on them as is the case on $B_3$). Moving to right (increasing $T$) along that segment, the tube size shrinks with $E$, while the trajectories progressively get closer to the tubes, eventually resulting in the origin of branch $B_3$. Other segments that connect the various branches $B_2$ through $B_{11}$ also behave similarly.

\section{Bifurcations in the MFG PDE, and comparison with the Hamiltonian BVP solutions}
\label{sec:bif_PDE}
In this section, we present some numerical solutions of the MFG Eqs. (\ref{eq:HJB1},\ref{eq:FP1}), and compare them with the solutions of the Hamiltonian BVP Eqs. (\ref{eq:Ham_ode}) discussed in the previous section, focusing on the $saddle\times center$ case. 

We note that the boundary conditions of the Hamiltonian BVP given by Eqs. (\ref{eq:ic_fc_ham}) prescribe the initial and final density, while in the standard MFG formulation, initial density and final value function are prescribed. We take the MFG initial and final conditions to match those used in the BVP, and use:
\begin{equation}\begin{split}
    \label{IC&FC: m}
    &m(x,0)=m_{IC}(x) \triangleq \frac{1}{\sqrt{2\pi \epsilon^2S^2(0)}}\exp\left(-\frac{(x-X(0))^2}{2\epsilon^2S^2(0)}\right),\\
    &m(x,T)=m_{FC}(x) \triangleq \frac{1}{\sqrt{2\pi \epsilon^2S^2(T)}}\exp\left(-\frac{(x-X(T))^2}{2\epsilon^2S^2(T)}\right),\end{split}
\end{equation} where $X(0)=-10,X(T)=10,S(0)=S(T)=4.5$ and $\epsilon=0.05$, as used in the Hamiltonian BVP problem.

The MFGs with prescribed initial and final density have been referred to as `MFG planning problems' \cite{porretta2014planning,chen2018steering} in the literature. The chosen final condition on density can be used to approximate the final value function, as described in \cite{achdou2012mean}. The idea is to employ the final value function condition $u(T,x)=\dfrac{1}{\epsilon_p}(m(x,T)-m_{FC}(x))$, for $\epsilon_p\ll 1$. With this approximation, the MFG PDE is reduced to a form which can be solved using different variants of Picard-Newton type algorithms \cite{lauriere2021numerical}. We employ one such variant to solve the problem at hand, and relegate its  details to the Appendix.

\begin{figure}
\centering
\begin{tabular}{c c}
\includegraphics[width=.49\textwidth]{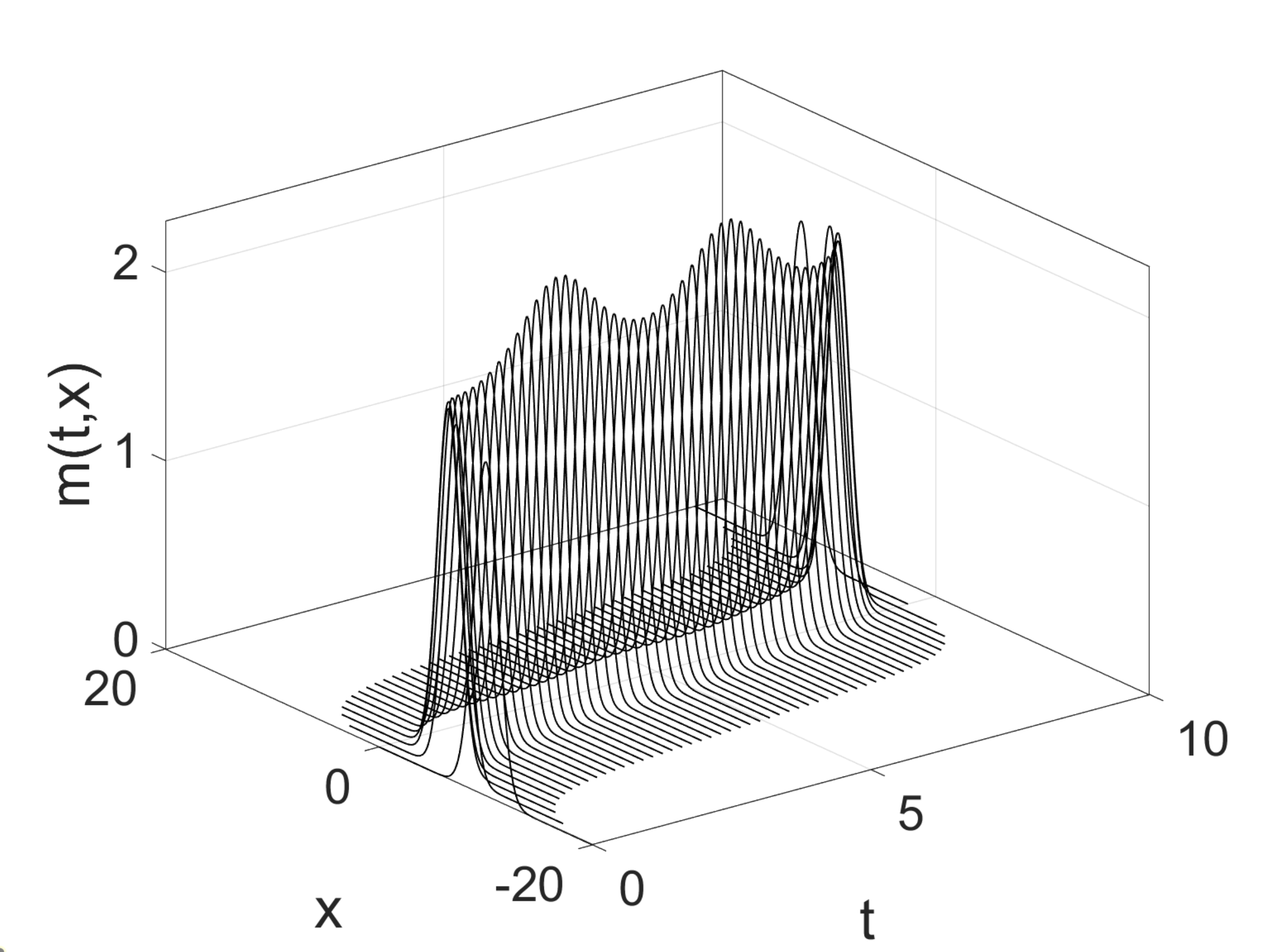} &\includegraphics[width=.49\textwidth]{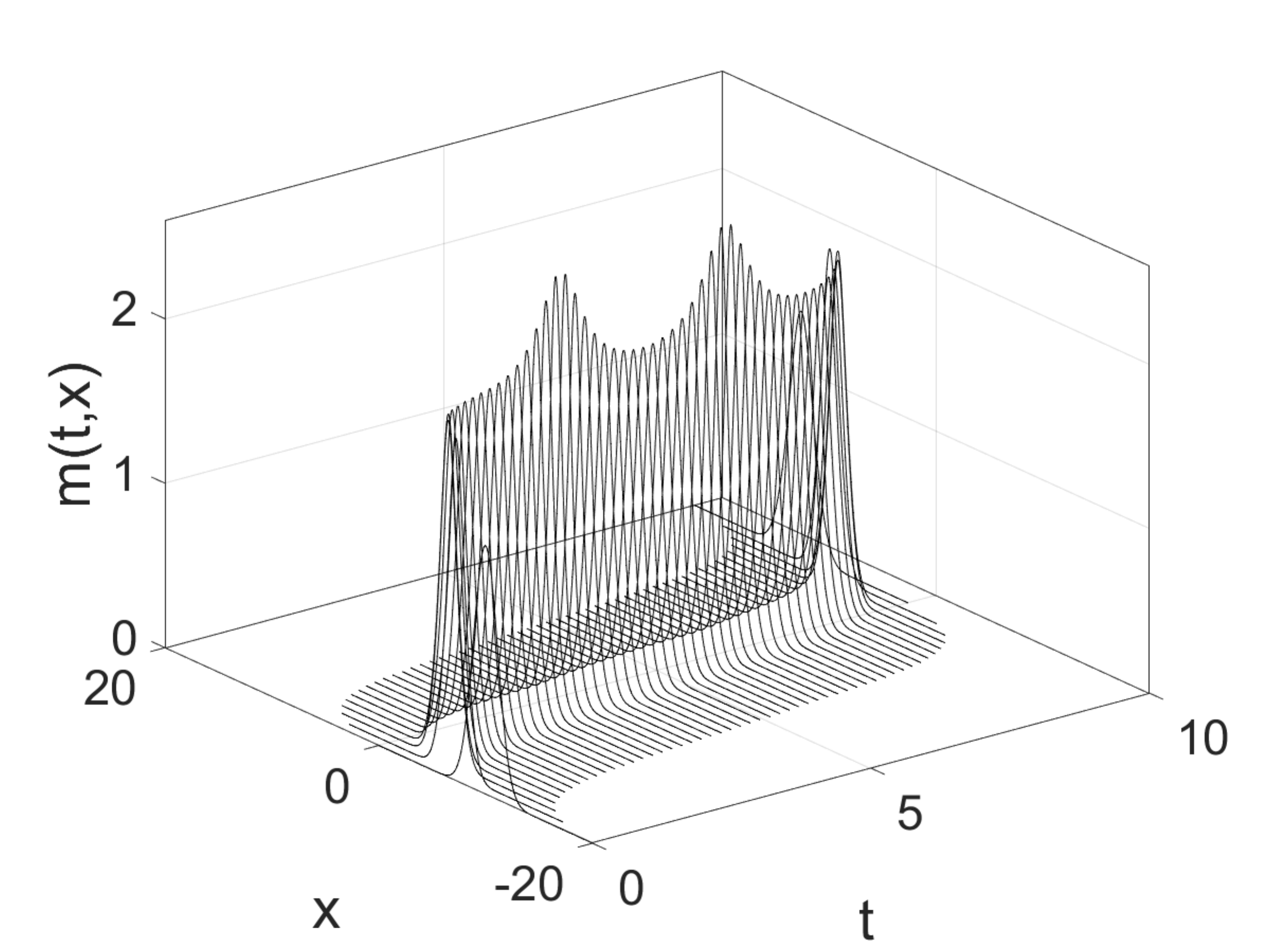} 
\\ (a) & (b)\\ 
\includegraphics[width=3in]{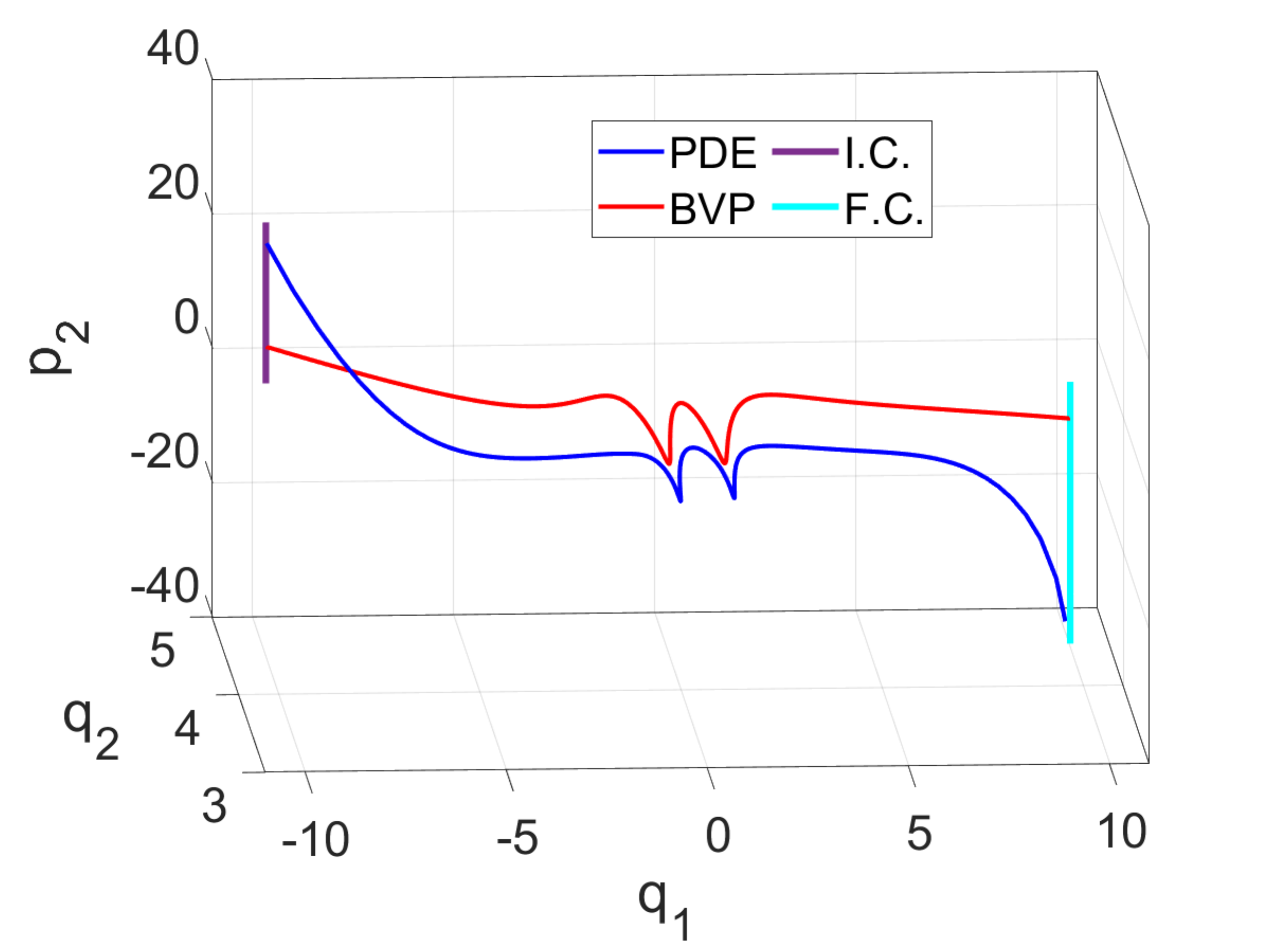}& \includegraphics[width=3in]{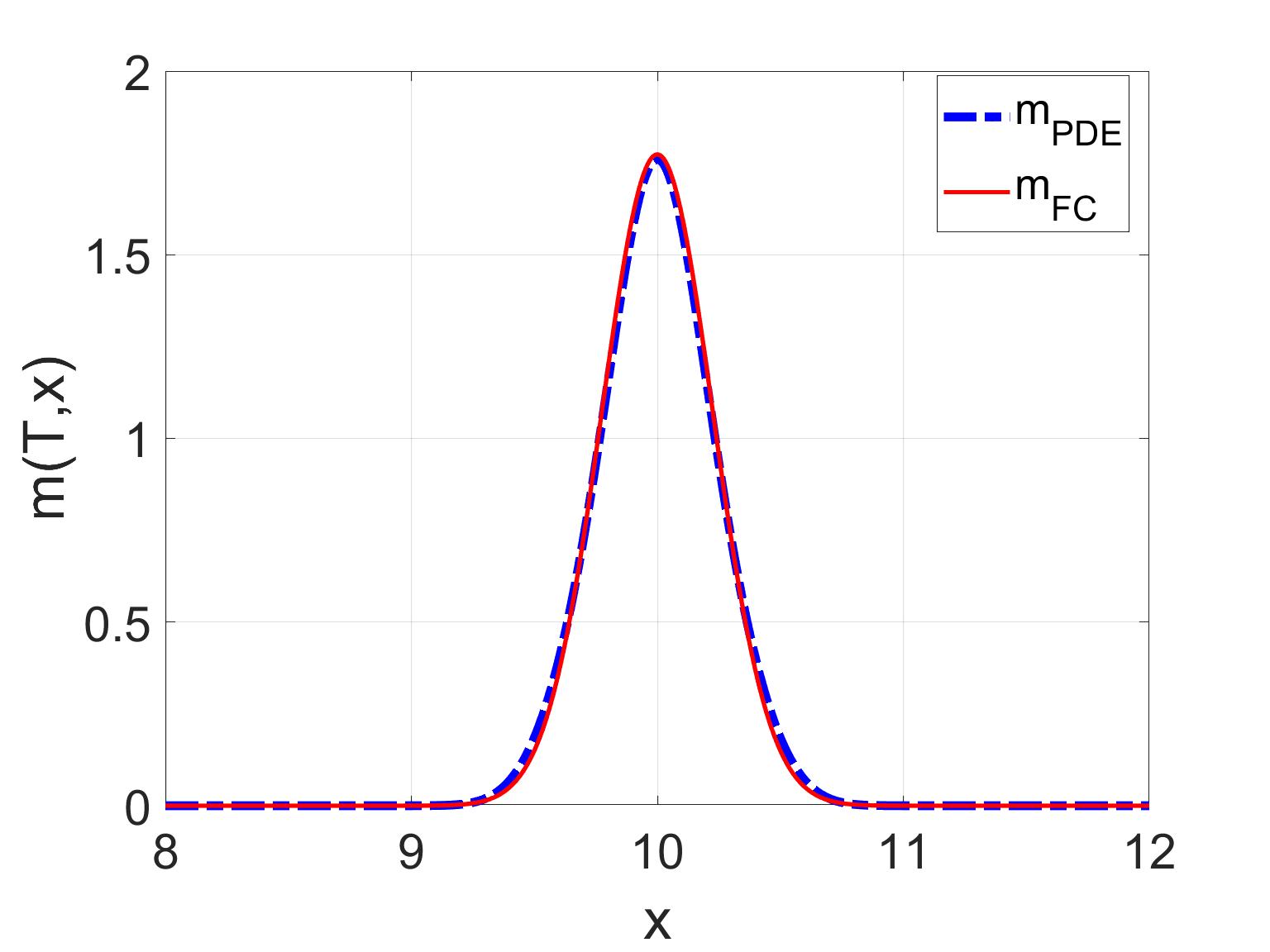} 
\\ (c) & (d)\\
\end{tabular}
\caption{\footnotesize The density $m(x,t)$ computed by solving the (a) Hamiltonian BVP, and (b) PDE equations for the same set of system parameters. (c) The BVP and PDE solutions in the phase space, showing that both trajectories complete two rotations in the ergodic phase. (d) The density at final time ($m_{PDE}$) computed using the PDE equations, and the prescribed final density ($m_{FC}$). }
\label{fig:comp_b2}
\end{figure}

\begin{figure}
\begin{tabular}{c c}
\includegraphics[width=0.48\textwidth]{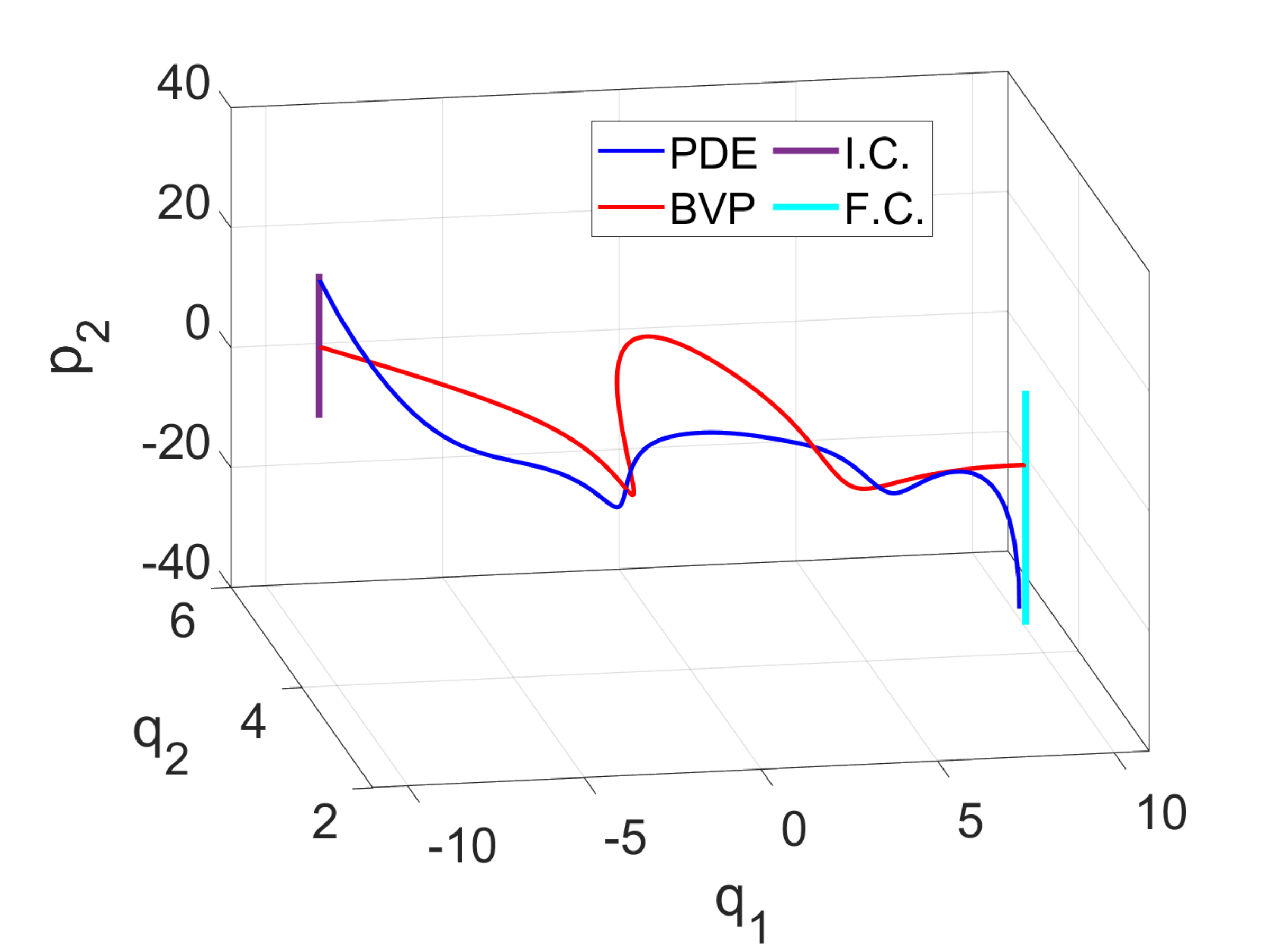} & \includegraphics[width=0.48\textwidth]{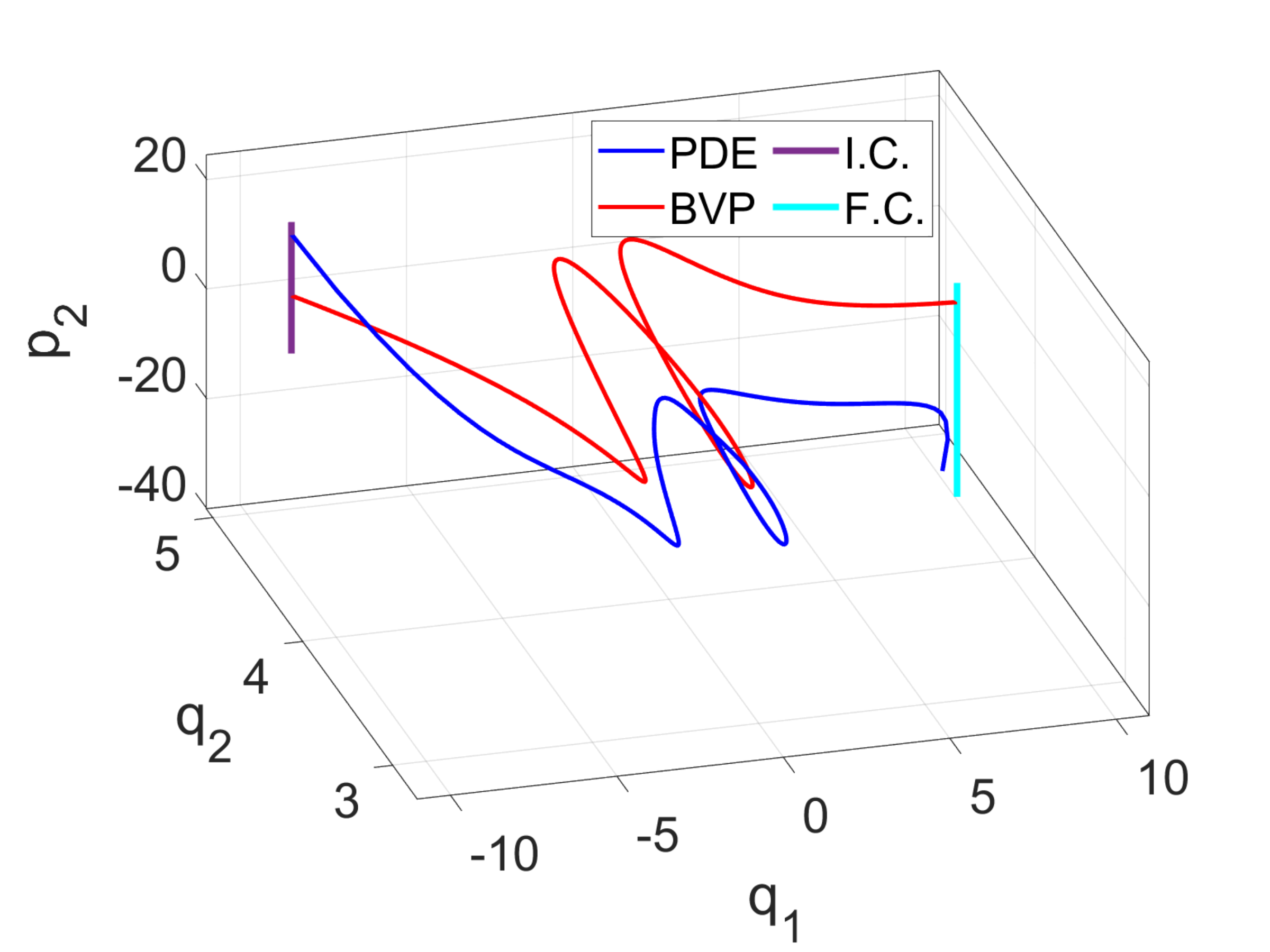} 
\\ (a) & (b)\\ 
\includegraphics[width=0.48\textwidth]{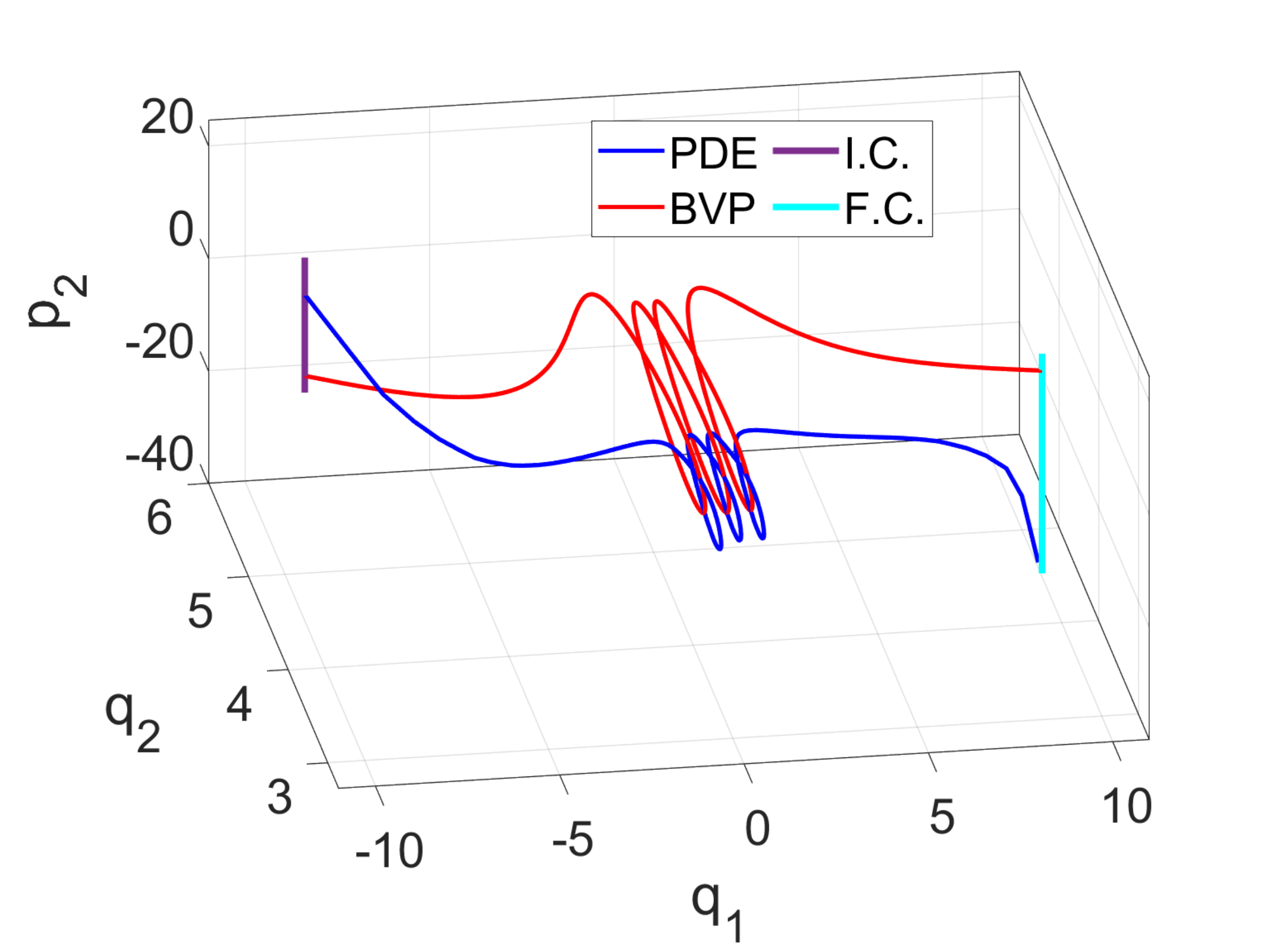} & 
\includegraphics[width=0.48\textwidth]{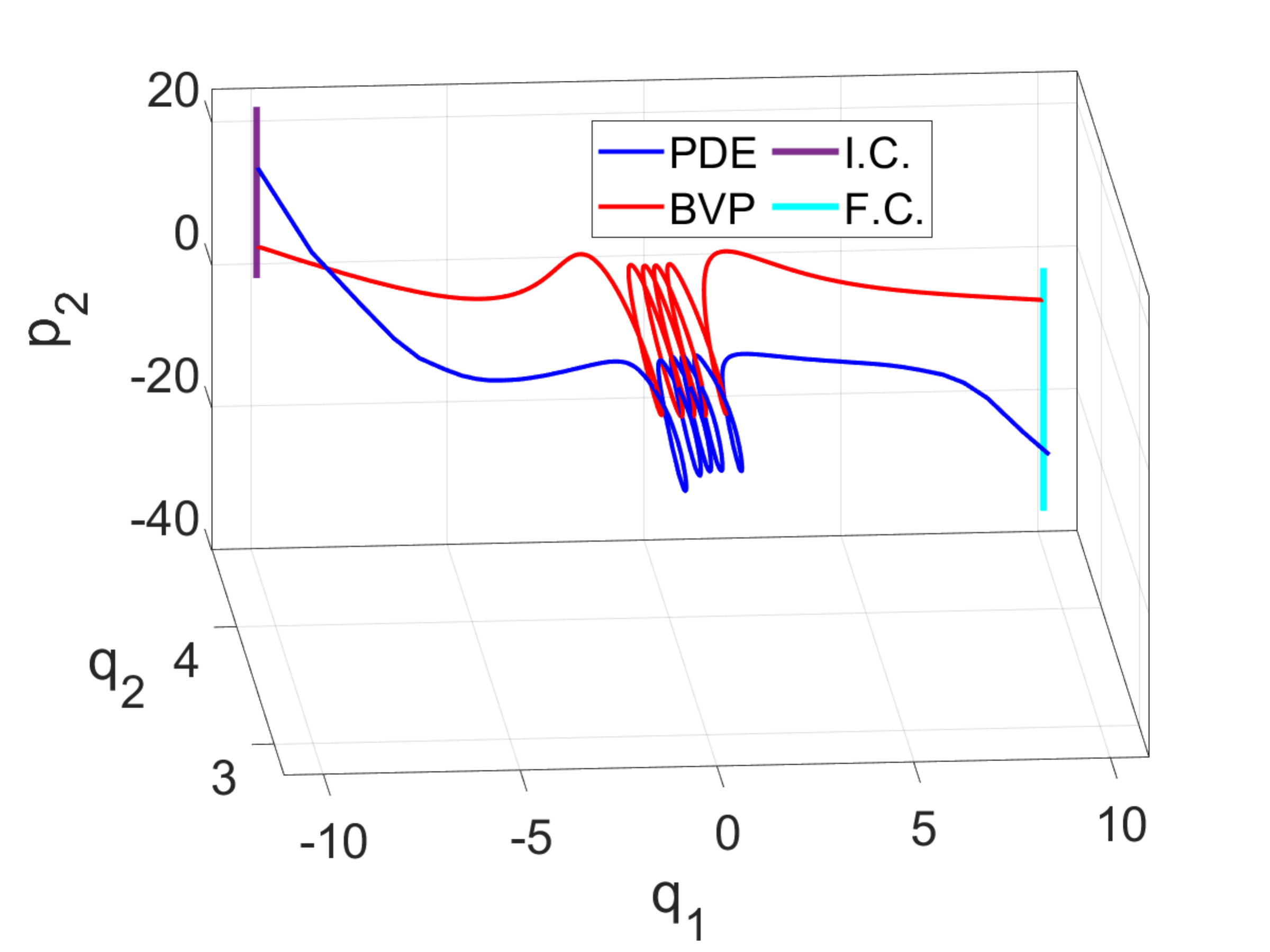} \\ (c) & (d)\\
\end{tabular}
\caption{\footnotesize Phase space portrait showing BVP and PDE solutions with (a) One, (b) one and a half, (c) three, and  (d) five rotations in the ergodic phase, respectively. Each pair of trajectories (BVP and PDE) are computed using the same set of system parameters. }
\label{fig:comp_ball}
\end{figure}

In order to solve the PDE system, we first pick system parameters (and initial guesses) corresponding to the $B_5$ branch of Fig. \ref{fig:bif_main}, in which case the BVP solution undergoes two full rotations in the ergodic phase. Figs. \ref{fig:comp_b2}(a,b) show the density evolution obtained from solving the BVP and PDE problems. In both cases, two peaks of the density during the evolution are evident, hinting at the topological similarity between the PDE and BVP solutions. To confirm this, we plot the two solutions in the phase space. To obtain phase space quantities from the PDE solution $(m(x,t),u(x,t))$, we employ the definitions of $(X,S,P,\Lambda)$ given in Sec. \ref{subsec:deriv_BVP}. Finally, we use the Legendre transformations of Sec. \ref{subsec:deriv_ham} to obtain  $(q_1(t),p_1(t),q_2(t),p_2(t))$. The phase space plots shown in Fig. \ref{fig:comp_b2}(c) confirm that the PDE and BVP solutions indeed have the same topology, i.e., they both undergo two full rotations in the ergodic phase. We observed similar agreement between the PDE and the BVP solutions for the other solutions branches, see Fig. \ref{fig:comp_ball}. Fig. \ref{fig:comp_b2}(d) shows that the PDE solution satisfies the final density condition.

\section{Discussion and Conclusions}
In this work, we have used a combination of reduced-order modeling and phase space analysis to identify and explain the topological nature of multiple solution branches of a variational non-monotonic finite-horizon MFG. This analysis rests on two key ideas: \begin{enumerate}
    \item The turnpike property of the MFG which allows the use of an ansatz borrowed from nonlinear Schrodinger equations to yield a 4D Hamiltonian BVP when the solution has solitonic form, and
    \item The organization of phase space in 4D Hamiltonian ODEs near a $saddle\times center$ equilibrium by invariant manifolds of the hyperbolic periodic orbit around that equilibrium point. 
\end{enumerate}
The phase space geometry near the equilibrium point of the Hamiltonian ODE is the key to determining the corresponding BVP trajectories. We first confirm that there is a unique branch of solutions in the $saddle\times saddle$ case with $E$ monotonically decreasing with increasing $T$, reminiscent of the uncoupled case treated in \cite{ullmo2019quadratic}. 

In the $saddle\times center$ case that is the focus of this work, the trajectories in the system linearized around the equilibrium can be classified as transit or non-transit, depending upon whether they are inside or outside the cylindrical invariant manifolds (tubes) of the periodic orbit, respectively. Due to hyperbolic nature of the periodic orbit, the conclusions from the linearized system persist in the nonlinear regime. This classification provides a characterization of initial and final conditions that can be joined via BVP trajectories for a given energy level $E$. In addition to a branch similar to that of the $saddle\times saddle$ case, several new branches, each with a fixed non-trivial topology, are identified. On these branches, the energy increases with $T$ while trajectories travel on the tubes. The topology is unambiguously determined by counting the number of intersections of the trajectory with a surface of section in the phase space, and equivalently, the number of half-rotations around the tube during the ergodic phase. 
The system parameters identified via BVP analysis are used while solving the full order MFG PDEs. The solutions obtained from the PDEs have the same topology as the corresponding BVP solutions, hence validating the model reduction methodology and the related phase space analysis.

This work adds to the small but growing corpus of results \cite{ullmo2019quadratic,hongler2020mean,bonnemain2020universal}, on the use of exact solutions, reduced-order modeling and dynamical systems analysis for understanding the behavior mean field games in various parametric regimes. It also demonstrates the role of geometry of invariant manifolds of periodic orbits in determining solutions of BVPs with a forward-backward nature. This could be of independent interest in problems related to optimal transport \cite{chen2021optimal,elamvazhuthi2018optimal1} and mean field control \cite{fornasier2014mean,carmona2013control}.

The results obtained here can be generalized in various directions. The tube dynamics based analysis described in this work can be generalized to (2n+2)D Hamiltonian ROMs with a rank-1 saddle (i.e., with a $saddle\underbrace{\times center \times center \times center}_{n\:times}$ type equilibrium point) for $n>1$. From the MFG viewpoint, this will allow the analysis of ROMs with more than two degrees of freedom, e.g., those with the controlled dynamics restricted to first $n+1$ moments of the distribution for $n>1$. 

Our framework can potentially also be used to construct homoclinic, heteroclinic, and `brake orbits' in certain classes of MFG systems \cite{cesaroni2021brake}. More generally, by considering reduced-order model regimes with multiple coexisting $saddle\times center$ equilibrium points, the configurations space ($q_1,q_2$) can be divided into various realms, with each pair of neighboring realms separated by a bottleneck around a $saddle\times center$ equilibrium point. Using prior results on symbolic dynamics in such systems \cite{koon2000heteroclinic}, one could construct BVP solutions with arbitrary itineraries, i.e., trajectories that visit different realms in any specified order. Some of these extensions will be taken up in future work.

\clearpage
\appendix*
\label{app}
\section{Numerical scheme for solving the MFG PDE system}

Consider the MFG system of Eqs.(\ref{eq:HJB1},\ref{eq:FP1}) :
\begin{equation}
    \begin{split}
        \label{mfg_p}
        &\partial_tu(x,t)-\frac{1}{2\mu}(\partial_xu(x,t))^2+\frac{\sigma^2}{2}\partial_{xx}u(x,t)=\bar{V}[m](x,t),
    \\ & \partial_tm(x,t)-\frac{1}{\mu}\partial_x(m(x,t)\partial_xu(x,t))-\frac{\sigma^2}{2}\partial_{xx}m(x,t)=0,
    \end{split}
\end{equation}
with prescribed initial density $m(x,0)=m_{IC}(x)$, and the final value function condition $u(T,x)=\dfrac{1}{\epsilon_p}(m(x,T)-m_{FC}(x))$.

We use finite-difference discretization in space and time to solve the above problem in the space-time domain $[-L/2,L/2]\times[0,T]$, where we pick $L$ to be large enough to avoid any boundary effects. The spatial interval is divided into $N_x$ uniform subintervals of size $\delta_x=L/N_x$, and the time interval $[0,T]$ is divided into $N_t$ time steps of size $\delta_t=T/N_t$. With this discretization, we use the notation $M^n_i\triangleq m(x_i,t_n)$, and $U^n_i\triangleq u(x_i,t_n)$, where $x_i$ is the $i$th spatial grid point, and $t_n$ is the $n$th temporal grid point. Also, $M^n\triangleq M^n_{i=0:N_x}$, and $U^n\triangleq U^n_{i=0:N_x}$. For all points other than those on the boundary of the spatial domain, we use central difference approximation for spatial derivatives. The forward difference scheme is used for time discretization. For $Z=M,$ or $Z=U$, finite difference operators are defined as follows:
\begin{equation}
    \begin{split}
        \label{D operators}
        &D_xZ^n_i=\frac{Z^n_{i+1}-Z^n_{i-1}}{2\delta_x},\\&        \Delta_xZ^n_i=\frac{Z^n_{i+1}-2Z_i^n+Z^n_{i-1}}{\delta_x^2},\\&
        D_tZ^n_i=\frac{Z^{n+1}_{i}-Z^n_{i}}{\delta_t}.
    \end{split}
\end{equation}

 Free boundary conditions are applied to the left ($i=0$) and right ($i=N_x$) boundary points, and we employ forward difference and backward difference scheme for spatial derivatives at the two locations, respectively.

The discretized versions of Eqs.(\ref{mfg_p}) are:
\begin{equation}
    \begin{split}\label{eq:MFG_dis}
        &-D_tU_i^n-\frac{\sigma^2}{2}\Delta_xU^n_{i} +\frac{1}{2\mu}(D_xU^n_{i})^2=-f(M^{n+1}_{i})-U_0(x_{i}),\\&
        D_tM_i^{n}-\frac{\sigma^2}{2}\Delta_xM^{n+1}_{i} -\frac{1}{\mu}\bigg((D_xM_i^{n+1})(D_xU_i^{n})+M^{n+1}_i\Delta_xU^n_i\bigg)=0,
    \end{split}
\end{equation}
with prescribed initial density $M^0=m_{IC}(x_{0:N_x})$, and final value function $U^{N_t}=\dfrac{1}{\epsilon_p}(M^{N_t}-m_{FC}(x_{0:N_x}))$. Following \cite{lauriere2021numerical}, the forward-backward discretized MFG system of Eqs.(\ref{eq:MFG_dis}) is solved using a Picard-type iteration. Given the density $\tilde{\mathcal{M}}^{[k]}\triangleq M^{0:N_t}_{0:N_x}$ at end of $k$th iteration, the value function $\mathcal{U}^{[k+1]}\triangleq U^{0:N_t}_{0:N_x}$ is obtained by solving the discretized nonlinear HJB repeatedly, starting at final time $T$, and marching backward in time up to $t=0$. At the $n$th time step, the system of equations to be  solved for $U^n$ is given by:
\begin{equation}\label{eq:FD_HJB}
    F_1(U_i^n)=-\frac{U_i^{n+1}-U_i^n}{\delta_t}-\frac{\sigma^2}{2}\Delta_xU^n_{i} +\frac{1}{2\mu}(D_xU^n_{i})^2+f(M^{n+1}_{i})+U_0(x_{i})=0,
\end{equation}
for $i=0:N_x$. We use the Newton-Raphson method to solve this system.

Using $\mathcal{U}^{[k+1]}$, the density $\mathcal{M}^{[k+1]}$ is obtained by solving the linear discretized FP equation repeatedly, starting at $t=0$, and marching forward in time to $t=T$. The linear system to be solved at the $n$th time step for $M^{n+1}$ is:
\begin{equation}\label{linearFP}
    A^{FP}M^{n+1}=M^{n}/\delta_t,\quad\qquad A_{q,i}^{FP}=\begin{cases} \dfrac{1}{\delta_t}+\dfrac{ \sigma^2}{\delta_x^2}-\dfrac{1}{\mu}\Delta_xU_i^n, & i=q   \\
                     -\dfrac{ \sigma^2}{ 2\delta_x^2}+\dfrac{1}{2\mu\delta_x}D_xU^n_{i}, &  i=q-1
                     \\ -\dfrac{\sigma^2}{ 2\delta_x^2}-\dfrac{1}{2\mu\delta_x}D_xU^n_{i}, &  i=q+1
                     \\ 0, & otherwise.
       \end{cases}.
\end{equation}
At the end of $(k+1)$th Picard iteration, we use damping to obtain the updated values of value function $\tilde{\mathcal{U}}^{[k+1]}$ and density $\tilde{\mathcal{M}}^{[k+1]}$. The Picard iterations are continued till tolerance is met, see Algorithm \ref{alg:Picard} for details.

\begin{table}[htbp]
\centering
\caption{\footnotesize Notation used in the pseudo-code for solving the MFG PDEs}
\begin{tabular}{|c|p{10cm}|}
\hline
\textbf{Notation} & \textbf{Description} \\
\hline
$k$ & Picard iteration index. \\
$k_{max}$ & Maximum number of Picard iterations allowed\\
$n$ & Time index. \\
$Tol$ & Tolerance for Picard iterations. \\
$N_x$ & Number of spatial subintervals. \\
$N_t$ & Number of time steps. \\
$\delta$ & Damping coefficient for Picard iteration. \\
$\delta_t$ & Time stepsize. \\
$\delta_x$ & Length of the each spatial subinterval. \\
$\tilde{\mathcal{M}}^{[k]}\in \mathbb{R}^{(N_t+1)\times(N_x+1) }$ & Density at end of $k$th Picard iteration. \\
$\tilde{\mathcal{U}}^{[k]}\in \mathbb{R}^{(N_t+1)\times(N_x+1) }$ & Value function at end of $k$th Picard iteration.\\
$\mathcal{U}^{[k+1]}\in \mathbb{R}^{(N_t+1)\times(N_x+1) }$ & Value function obtained  from solving the HJB during $(k+1)$th Picard iteration. \\
$\mathcal{M}^{[k+1]}\in \mathbb{R}^{(N_t+1)\times(N_x+1) }$ & Density obtained from solving the FP during $(k+1)$th Picard iteration. \\

\hline
\end{tabular}
\label{tab:notation_mfg}
\end{table}

\begin{algorithm}
\caption{\footnotesize Pseudo-code for solving the MFG PDEs}\label{alg:Picard}
\begin{algorithmic}
    \State Make initial guesses for the density $\tilde{\mathcal{M}}^{[0]}$ and the value function $\tilde{\mathcal{U}}^{[0]}$.
    \State $isconverged=0$
    \For{$k = 0$ to $k_{max}$} \Comment{Picard iteration}
        \State  $M^{0:N_t}=\Tilde{\mathcal{M}}^{[k]}$
        \State Set final time value function $U^{N_t}=\dfrac{1}{\epsilon_p}(M^{N_t}-m_{FC}(x_{0:N_x}))$
        \For{$n = N_t-1$ down to $0$} \Comment{HJB time marching}
            \State Solve Eq. \ref{eq:FD_HJB} using Newton-Raphson to get ${U}^{n}$
        \EndFor
        \State $\mathcal{U}^{[k+1]}= U^{0:N_t}$
      \State Set initial time density $M^0=m_{IC}(x_{0:N_x})$
        \For{$n = 0$ to $N_t-1$}            \Comment{FP time marching}
            \State Solve the linear system of Eq. (\ref{linearFP}) to get $M^{n+1}$
        \EndFor
         \State $\mathcal{M}^{[k+1]}= M^{0:N_t}$
        \State Update:
        \State $\tilde{\mathcal{M}} ^{[k+1]} = \delta^{(k)}\tilde{\mathcal{M}} ^{[k]}+(1-\delta^{(k)})\mathcal{M}^{[k+1]}$
        \State $\tilde{\mathcal{U}} ^{[k+1]} = \delta^{(k)}\tilde{\mathcal{U}}^{[k]}+(1-\delta^{(k)})U^{(k+1)}$
        \If{$||\tilde{\mathcal{U}} ^{[k+1])} - \Tilde{\mathcal{U}}^{[k]}|| < Tol$ And $||\Tilde{\mathcal{M}}^{[k+1]} - \Tilde{\mathcal{M}}^{[k]}|| < Tol$}
        \State $isconverged=1$
        \State break
        \EndIf
    \EndFor
\end{algorithmic}
\end{algorithm}

\subsection{Convergence}

\begin{figure}[hbt!]
\centering
\includegraphics*[width=.75\textwidth]{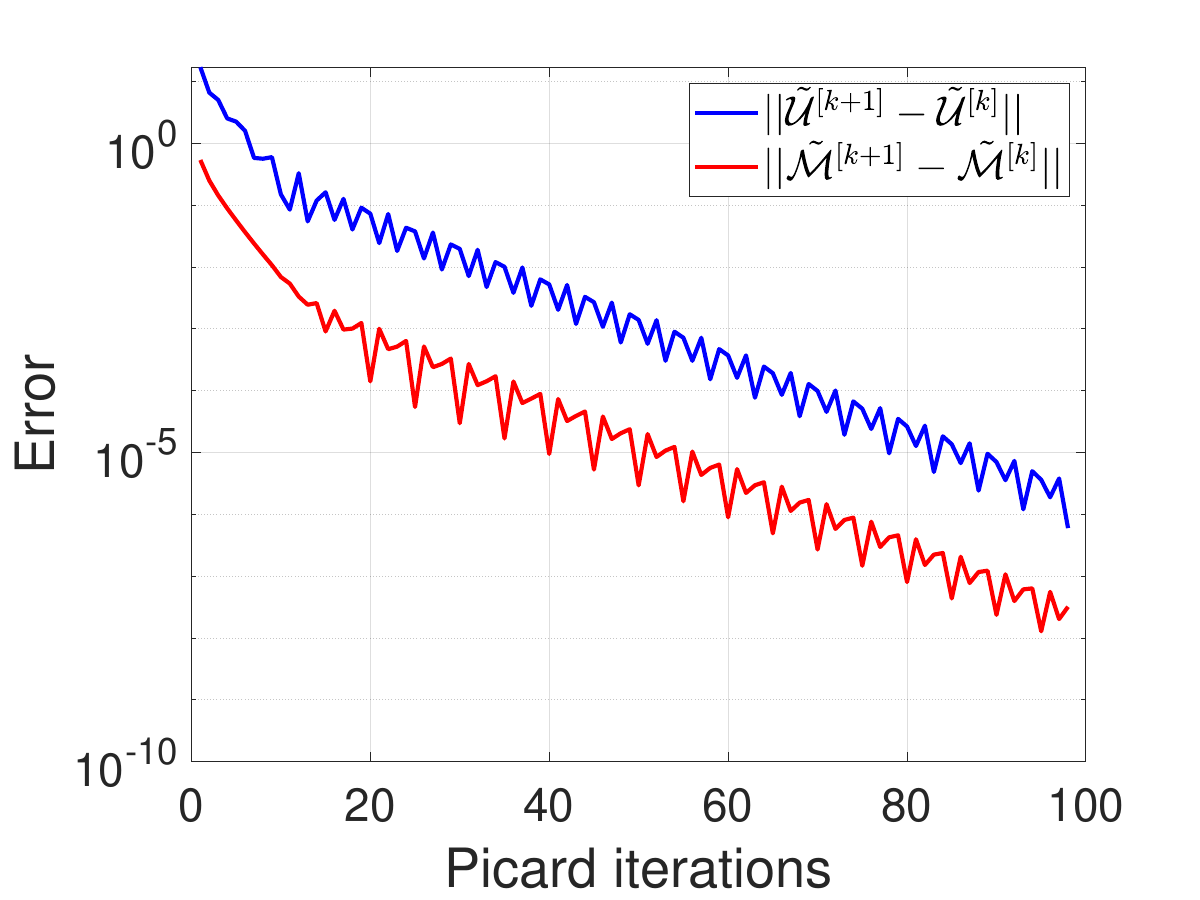} 
\caption{\footnotesize Decay of the two error terms in the Picard iteration Algorithm \ref{alg:Picard} to solve the MFG PDEs for the case shown in Fig. \ref{fig:comp_b2}. } 
\label{fig:picard_conv}
\end{figure}

For solving the MFG PDEs for the two rotation case discussed in Sec. \ref{sec:bif_PDE}, we used the following parameter values for the algorithm: $L=40$, $N_x=500$, $N_t=500$, $T=9.5$, $\delta=0.5$, $k_{max}=1000$, $\epsilon_p=0.01$, and $Tol=10^{-6}$. This yields $\delta_t=0.019$ and $\delta_x=0.08$. The convergence behavior of Algorithm \ref{alg:Picard} for this case is shown in Fig. \ref{fig:picard_conv}. It is evident that the algorithm converges rapidly. Similar performance was observed for solutions with different topologies (e.g., those shown in Fig. \ref{fig:comp_ball}).

\clearpage
\section*{Acknowledgments}
This material is based upon work supported by the US National Science Foundation under Grant No. 2102112.

\clearpage

\bibliographystyle{unsrt}
 \bibliography{refs}

\end{document}